\documentclass{nm}
\renewcommand\thisnumber{x}
\renewcommand\thisyear {201x}

\renewcommand\thisvolume{xx}
\usepackage{charter}
\usepackage[charter]{mathdesign}
\usepackage{amsmath}
\usepackage{graphicx}%
\usepackage{latexsym ,rawfonts}
 \usepackage{multirow}
\usepackage{tikz}
\usetikzlibrary{shapes}
\usepackage{color}

\usepackage{algorithmic}
\usepackage{algorithm}
\usepackage{subfigure}
\newcommand{\bi}{\begin{itemize}}
\newcommand{\ei}{\end{itemize}}

\newcommand{\ta}{\tilde{A}}
\newcommand{\caliber}{{c}}
\newcommand{\xk}{x^{(k)}}
\newcommand{\finevar}{F}
\newcommand{\coarsevar}{C}
\newcommand{\tvV}{\mathcal{V}}

\newcommand{\innerprod}[3][2]{\langle #2,#3 \rangle_{#1}}

\begin{document}

\markboth{A. Brandt, J. Brannick, K. Kahl, and I. Livshits}{Bootstrap Algebraic Multigrid}
\title{Bootstrap Algebraic Multigrid: status report, open problems, and outlook}
\author[Achi Brandt, James Brannick, Karsten Kahl, and Ira Livshits]{Achi Brandt\affil{1}, James Brannick\affil{2}\comma\corrauth, Karsten Kahl\affil{3}, and Ira Livshits\affil{4}}
\address{\affilnum{1}\ Microsoft Corporation, Mountainview, CA, 94043, USA.
 \\
 \affilnum{2}\ Department of Mathematics, The
       Pennsylvania State University, University Park, PA 16802, USA.
\\
 \affilnum{3}\ Department of Mathematics, University of Wuppertal, Wuppertal, 
 Germany.
 \\
  \affilnum{4}\ Department of Mathematical Sciences, Ball State University, Muncie IN, 47306, USA.
       }

\emails{
 {\tt acbrandt@microsoft.com} (Achi Brandt),
  {\tt brannick@psu.edu} (James Brannick),
 {\tt kkahl@math.uni-wuppertal.de} (Karsten Kahl),
 {\tt ilivshits@bsu.edu} (lra Livshits)
 }

\begin{abstract}
This paper provides an overview of the main ideas driving the bootstrap algebraic multigrid methodology, including compatible relaxation and algebraic distances for defining effective coarsening strategies, the least squares method for computing accurate prolongation operators and the bootstrap cycles for computing the test vectors that are used in the least squares process.
We review some recent research in the development, analysis and application of bootstrap algebraic multigrid and point to open problems in these areas.  Results from our previous research as well as some new results for some model diffusion problems with highly oscillatory diffusion coefficient are presented
to illustrate the basic components of the BAMG algorithm.  
\end{abstract}




\maketitle
\section{Introduction}
{Numerous scientific and engineering problems find their formulation in terms of (systems)
of partial differential equations, which in turn require the solution of 
large-scale finite element, finite difference, or finite volume equations. 
Modern applications involve large-scale parallel processing of
(linear) systems with millions or even billions of unknowns, for which multigrid (MG) 
methods often provide solvers that are optimal with respect to their computational complexity
and, hence, their parallel scalability.  

Multigrid solvers for sparse systems of linear equations
$Ax=b$ are all based on two complimentary processes:
a local relaxation scheme (smoother) that eliminates certain
components of the error by working on the fine level and a coarse-level 
correction that treats the remaining global error.} 

Generally, the design and analysis of these two MG components are 
based on the following smoothing property of relaxation.
For any $0 < \rho < 1$, an error vector $e$ is called $\rho$-smooth if
all its {\it normalized} residuals\footnote{The vector of {\it
normalized residuals} is $\ta e$, where $\ta$, the {\it normalized
matrix}, is scaled so that the $l_1$ norm of each of its rows is
$1$. Correspondingly, we define the {\it normalized eigenvalues}
of $A$ to be the eigenvalues of $\ta$. The magnitude of any
normalized eigenvalue is at most $1$.} are smaller than $\rho
\|e\|$. The basic observation \cite{B83} is that the convergence
of a proper relaxation process\footnote{Kaczmarz relaxation can
always be used (but better schemes are often available),
supplemented when needed by local relaxation steps around
exceptionally large residuals (as suggested in \cite[\S
A.9]{BAMG1}).} slows down only when the current error is
$\rho$-smooth with $\rho \ll 1$, the smaller the $\rho$ the slower
the convergence. This basic observation implies that when
relaxation slows down, the error vector $e$ can be approximated in
a much lower-dimensional subspace. Very efficient {\it ``geometric
multigrid''} solvers have been developed for the case that the
lower-dimensional subspace corresponds to functions on a
well-structured grid (the {\it coarse} level), on which the smooth
errors can be approximated by easy-to-derive equations, based for
example on discretizing the same continuum operator that has given
rise to the fine-level equations $Ax=b$ {to define the coarse-level operator
$A_c$.} The coarse-level
equations are solved using recursively a similar combination of
relaxation sweeps and still-coarser-level approximations to the
resulting smooth errors.

The basic two-grid method for solving $Ax=b$, from which a multigrid method is defined by recursion,
involves a stationary linear iterative 
method (the smoother) applied to the fine-grid system,
and a coarse-grid correction: given an approximation $w\in \mathbb
C^n$, compute an update $v\in \mathbb C^n$ by \medskip
\begin{enumerate}
\item Pre-smoothing: $y = w+M(f-Aw)$,
\item Correction: $v = y + PA_c^{-1}R(f-Ay)$.
\end{enumerate}\medskip
Here, $M$ is the approximate inverse of $A$ that defines the multigrid smoother and
$R: \mathbb{C}^{n} \mapsto \mathbb{C}^{n_c}$ and $P: \mathbb{C}^{n_c} \mapsto \mathbb{C}^n$ with $n_c < n$ are the restriction and 
interpolation operators that map information between the coarse grid of size $n_c$ and the
fine grid of size $n$. 

To deal with much more general situations, where the
fine-level system may not be defined on a well-structured grid
nor perhaps arise from any continuum problem, {\it ``algebraic
multigrid''} (AMG) methods were developed to derive the
set of coarse-level {\it variables} and coarse-level {\it
equations} directly from the given matrix $A$. The basic
approach (developed in \cite{BMR83, BMR84, RS87} and called today
``classical AMG'' or RS-AMG) involves the following two steps. \bi
\item[(1)] Choosing the coarse-grid variables, e.g., as a subset $C$ of
the set of fine-grid variables, $\Omega$, such that each variable in $F=\Omega
\setminus C$ is
``strongly connected'' to variables in $C$. \item[(2)]
Approximating the fine-level residual problem $Ae=r$ by the
coarse-level equations $A_c e_c=r_c$ using the Galerkin
definition $A_c = RAP$ with $r_c=Rr$, so that $e \approx P e_c$.\footnote{{
A common choice when $A$ is HPD is to set $R=P^{H}$, which leads to a variational
Galerkin scheme.}}
\ei 
In this classical approach, the ``strength of connection'' measure as well as the
interpolation matrix $P$, and the {restriction matrix} $R$
are all defined directly in terms of the elements of the matrix
$A$.
The overall algorithm gives satisfactory results only when the (properly scaled) matrix $A$ has a dominant diagonal and (with small possible exceptions) all its off-diagonal elements have the same sign. Even then, the produced interpolation in many cases has limited accuracy, insufficient for full multigrid efficiency.

{The Bootstrap AMG (BAMG) algorithm, introduced in \cite[\S
17.2]{B00}, and the closely related Adaptive AMG algorithms introduced in~\cite{MBrezina_RFalgout_SMacLachlan_TManteuffel_SMcCormick_JRuge_2003,MBrezina_2005}, were developed 
to extend AMG to an even wider class of problems for which the assumptions on smooth error
made in classical AMG are violated.
The main idea in Bootstrap and Adaptive AMG is to 
apply the current solver to appropriately formulated
test problems to explicitly compute approximations of the
smooth error components.
These components are then incorporated
into the coarse space by the AMG setup algorithm 
as needed.   

The adaptive algebraic multigrid methods developed originally in \cite{MBrezina_RFalgout_SMacLachlan_TManteuffel_SMcCormick_JRuge_2003,MBrezina_2005},
use a sequential setup procedure to adaptively identify the components of the space that need to be represented in the range of the coarse grid operator.  While such an approach insures that test vectors are incorporated into the coarse space only as they are needed, this filtering procedure requires reconstructing the entire multilevel hierarchy at each iteration, i.e., recomputing interpolation and the coarse grid matrices on all levels.  Thus, these approaches are efficient when one (or few) components are needed to define an effective interpolation operator, but can be costly in cases where many such components are needed to obtain full MG efficiency.

The two main ideas in bootstrap AMG that distinguish the approach from other Adaptive AMG methods are: (1) the idea of defining restriction and interpolation to fit a collection of smooth vectors in a least squares sense, and (2) the idea of computing the collection of test vectors simultaneously in a multigrid scheme, where the bulk of the work is done on coarse, i.e., on cheap levels. 
The BAMG algorithm, described below, combines the general devices of (A) compatible relaxation and algebraic distances to choose the coarse variables
with (B) a least squares procedure and (C) a bootstrap cycling scheme for computing the test vectors
used in the LS process to construct high-quality restriction and interpolation operators
and thereby a multigrid hierarchy.   
Before giving an overview of those Bootstrap techniques, we remark here that they should be used as needed, e.g., in cases where geometric MG and AMG solvers are not available and/or the additional costs they introduce can be amortized, such as problems where the same system has to be solver for multiple right hand sides.}

\bi
\item[(A)]  The coarse variables are selected using compatible relaxation
and algebraic distances, as tools for gauging the quality of the coarse variable set
and guiding the coarsening process, respectively.  \ei

A general criterion for measuring the quality of the set of coarse
variables $C$ is the fast convergence of compatible relaxation
(CR), as introduced in \cite{etna2000} (being a special case of
choosing coarse variables for much more general types of systems
\cite{SU}, introduced first for problems in statistical mechanics
\cite{RMG}).  Generally, CR refers to a relaxation scheme that keeps
the coarse-grid variables invariant.  {A general approach introduced by 
Brandt and Livne (see also~\cite{HCR}), that can be applied
to any relaxation scheme is 
``habituated CR'' (HCR), which uses the full fine-grid relaxation scheme and
then imposes certain constraints on the coarse-grid variables.
The HCR criterion allows choosing a set of coarse variables which
is not necessarily a subset of the set of fine-level variables. In
particular, coarse variables can represent some {\it averages} of
fine variables, as needed for example for equation homogenization
or for coarsening nonlinear problems. Another important type of a
set of coarse variables whose choice can be guided by CR arises in
{\it wave equations}, where it should include several different
subsets, each for example corresponding to a different wave
direction.
Moreover, given a suitable definition of the averages used in defining 
coarse variables, 
the HCR convergence rate yields an
accurate prediction of the convergence rate that can be achieved
by a two-grid solver for a given relaxation scheme and the
proposed set of coarse variables~\cite{JBrannick_2005a,JBrannick_RFalgout}. In this way, the HCR quantitative prediction
is also very useful in designing and debugging the actual solver.}

Algebraic distance gives a general approach 
for measuring the connectivity among nodes in a graph. 
In AMG, the algebraic distance measure
 is based on a notion of strength of connectivity among variables that is derived by minimizing a local least squares (LS) problem which is also used to derive interpolation (described in Section~\ref{sec:LS}). The approach first constructs direct (caliber-one) LS interpolation for a given set of test vectors (representatives of algebraically smooth error obtained by the bootstrap process described
in Section~\ref{sec:BS}) and then defines the algebraic distance between a fine point and its neighboring points in terms of the values of the local least squares functionals resulting from the so defined interpolation. The algebraic distance measure thus aims to address the issue of strength of connections in a general context -- the goal being to determine explicitly those degrees of freedom from which high quality least squares interpolation for some given set of smooth test vectors can be constructed.
  
An algorithm using algebraic distance was first developed for combinatorial optimization problems on graphs~\cite{SRB09,RonSB11}.   In these works, algebraic distance is developed as a measure of the closeness between nodes that is then used in an aggregation algorithm to build multilevel representations of graphs.
In~\cite{BAMG_Aniso}, algebraic distance was developed as a measure of strength of connectivity 
among variables for coarsening non-grid-aligned anisotropic diffusion problems.  Further, it was shown that combing the measure with compatible relaxation results in a robust algorithm for choosing 
interpolation operators for such problems.  And an affinity measure based on algebraic distance
together with an energy measure 
are used to form aggregates in the Lean AMG solver developed in~\cite{LAMG_Report} for the graph Laplacian problem.
Related work is found in~\cite{brannick_local_stab_2011}.

\bi
\item[(B)] {The interpolation operators (as well as restriction operators for non-Hermitian problems) are derived by a least-squares
fit to a set of $\rho$-smooth Test Vectors (TVs) with
sufficiently small $\rho$}. \ei 

Denote by $C_i$ the set of coarse variables used in
interpolating to a fine grid variable $i\in\finevar$. 
Then the main idea in the least squares interpolation approach is to compute for each
$i\in \finevar$ a row of interpolation that approximates
the given set of test vectors at the points in $C_i$, minimizing the interpolation error for these
vectors in a least squares sense.
It follows from the
satisfaction of the CR criterion that, with a proper choice of
$C_i$, there exists an interpolation that will have only $O(\rho)$
errors in reproducing $\xk_i$ for {\it any} normalized vector $x^{(k)}$ which
is $\rho$-smooth. The size $|C_i|$ of this set should in principle
increase as $\rho$ decreases, but in practice a preset 
interpolation {\it ''caliber''} $\caliber$ yields small enough
errors. For example, for a discretized scalar PDE in
$d$-dimensions, a proper set with $|C_i|=d+1$ points would give
$O(h^2+\rho)$ accuracy, which may suffice if $\rho$ is small
enough (see below). Typically, the number of TVs need not be larger than $1.5
\caliber$ for this accuracy to be attained for almost every
$\xk_i$. 

The set $C_i$ can often be adequately chosen by geometric
considerations, such as the set of geometrical neighbors. 
If the chosen set is inadequate, the
least-squares procedure will show {\it bad fitness} (interpolation
errors large compared with $\rho$), and the set can then be
improved. The least-squares procedure can also detect variables in
$C_i$ that can be discarded without a significant loss in accuracy.
Thus, this approach allows creating interpolation with whatever
needed accuracy which is {\it as sparse as possible}. Moreover, the LS
derivation of interpolation is applicable to {\it general types}
of coarse variables (e.g., averages). \bi \item[(C)] {The TVs
are constructed in a bootstrap manner}, in which several tentative
AMG levels generated by fitting interpolation to
moderately smooth TVs are used to produce {\it improved} (much
smoother) TVs, repeating this if needed until fully efficient AMG
levels have been generated. \ei

Quantitatively, the algorithm is driven by the following
relations.
The first TVs are each produced by relaxing the homogeneous system
$Ax=0$ with a different starting vector. This quickly leads to a
set of $\rho$-smooth TVs with $\rho \ll 1$.  A combination of random
vectors and additional geometrically derived smooth vectors can generally be
used as starting vectors.  In the case of discretized isotropic
PDEs, if geometrically smooth vectors that satisfy the homogeneous
boundary conditions are used, relaxation may not be needed at all.
In many other cases relaxation can be confined to the neighborhood
of boundaries and discontinuities. As mentioned, the
least-squares fitting to such a $\rho$-smooth set produces an
interpolation scheme whose errors in interpolating {\it any} normalized
$\rho$-smooth vector are at most $O(\rho)$. The coarse-level
Galerkin system based on such an interpolation will then yield an
$O(\rho^2)$-accurate approximation to every normalized eigenvalue
of a $\rho$-smooth eigenvector. Increasingly coarser AMG levels
can recursively be built in this way, all typically having such
$O(\rho^2)$ accuracy for eigenvectors with small normalized eigenvalues
and their linear combinations.

The AMG cycles based on these levels can further be applied to the
homogeneous system $Ax=0$. 
In many cases, the convergence (to
$x=0$) will exhibit the desired multigrid efficiency (close to the
HCR prediction), indicating that the produced AMG structure is
ready for production. If the given method needs further improvement, these cycles will quickly produce
$O(\rho^2)$-smooth TVs, which can then be used for a second round
of constructing the AMG hierarchy, potentially (with suitable
sets $C_i)$ providing $O(\rho^4)$-accuracy. More such rounds of TV creation, to obtain
still higher accuracy, is seldom needed.  

Some important
applications do result in nearly singular systems, e.g., the discretized Dirac equation arising in 
lattice quantum chromodynamics.
For such problems the system matrix has several almost 
zero modes (AZMs), which are a few special, nearly-singular eigenvectors 
with normalized eigenvalues much smaller than usual. For example, in discretizing a $q$-order PDE on a (structured or unstructured) grid with $O(h)$ meshsizes, the lowest normalized eigenvalues will all usually have size $O(h^q)$; but by a special adjustment of a parameter in the PDE, one normalized eigenvalue can be made arbitrarily small, having size $\eta \ll h^q$.

Usually, the potential existence of one, or sometimes several,
AZMs is known in advance, as a property of the problem at hand. To
obtain $O(\eta)$ approximations to normalized eigenvectors would
generally require substantially increasing the number of interpolation points, or caliber
denoted by $\caliber$. However, since the number of AZMs is typically small, the
interpolation can be tailored to fit them individually
sufficiently well, without increasing its caliber. This will in
fact automatically happen as soon as the AZMs emerge (in the
cycles applied to the homogeneous system), provided the
least-squares fitting weights each test vector 
proportionally to the inverse of its normalized residual norm.
We note that a multilevel eigensolver can be integrated into
the bootstrap setup process to further accelerate the convergence 
of the setup for the AZMs, as introduced in~\cite{BAMG2010}.  

In the remainder of the paper, an overview of the main components of the BAMG approach are given. 
Section \ref{sec:CR} contains a review of compatible relaxation and algebraic distances.  
The least squares approach for computing interpolation is discussed in Section \ref{sec:LS}
and the bootstrap procedure aimed at computing increasingly accurate TVs is
discussed in Section \ref{sec:BS}.  A main focus in the presentation in these 
sections is to point to open research questions in the development, analysis and
application of BAMG and to the most pressing research directions on these topics.
Promising future research areas in BAMG are outlined in Section \ref{sec:FD}.

\section{Compatible relaxation and algebraic distance}\label{sec:CR}
The first task in coarsening a given system is the choice of coarse-level variables, which in the BAMG
framework is accomplished using compatible relaxation and algebraic distances. 
Compatible relaxation coarsening algorithms have been developed and analyzed for both classical
AMG~\cite{HCR,PanayotRob_2003,JBrannick_RFalgout,Brannick_Trace_06} and aggregation-based coarsening~\cite{Jacob}.  These works have demonstrated the utility of
CR as a tool for guiding the coarsening process and assessing the quality of a given coarse variable set.  
The works in~\cite{HCR,JBrannick_RFalgout} focused on adding robustness to classical AMG using CR coarsening approaches, in particular for anisotropic diffusion problems with highly oscillatory and strongly anisotropic diffusion coefficient.  The method developed in~\cite{Jacob} does not use CR directly in the setup algorithm, however, CR was used as the main tool for the analysis and design of the resulting algorithm, which 
was developed specifically for general non-grid-aligned anisotropic diffusion problems.  

An outline of a generic coarsening process based on CR for the linear system $Ax=b$ is given in Algorithm 
\ref{alg:CR}.  Here, $\rho_{cr}$ is an approximation to the asymptotic convergence rate of CR computed after the $\nu$ iterations are applied and $\delta$ is a tolerance that can be set by the user or more generally it can be optimized within the algorithm.  In~\cite{HCR}, the CR coarsening algorithm is designed to approximately optimize $\beta = \rho_{cr}^{\frac{1}{W}}$, where $W$ denotes the total work of a single AMG cycle that is estimated using the  coarsening ratios in different stages of a CR coarsening algorithm.  
Starting with an initial set $C_0$ of coarse variables, whose values are held fixed at zero, a few sweeps of compatible relaxation (with vanishing right-hand side so that the solution should converge to zero, but starting of course from a suitable non-zero first approximation) are applied to detect convergence slowness if the coarse-variable set is inadequate. In this case, a set of variables that should be added to the coarse set needs to be
determined.   The typical approach is to select a subset of the set of variables which are slow to converge (to zero) and add these to the set of coarse variables.  As in AMG, the additional variables to be added to $C_0$ can be chosen using a notion of connectivity of a variable to its neighbors, which is defined via algebraic distance in BAMG. 

\begin{algorithm}
\caption{: Compatible relaxation coarsening, \textit{Input:} $\coarsevar_{0}$ \textit{Output:} $\coarsevar$\label{alg:CR}}
\begin{algorithmic} 
\STATE Initialize $\coarsevar = \coarsevar_{0}$
\STATE Perform $\nu$ CR iterations
\WHILE{$\rho_{cr} > \delta$}
\STATE Update $\coarsevar$ 
\STATE Perform $\nu$ CR iterations 
\ENDWHILE
\end{algorithmic}
\end{algorithm}

{In order to illustrate the flow of the above CR coarsening process, we report results from~\cite{JBrannick_RFalgout} for a specific implementation of the algorithm applied to the anisotropic diffusion problem
\begin{equation}\label{eq:jumps}
 - \nabla \cdot \mathcal{K}(x,y) \nabla u(x,y) = f(x,y),
 \end{equation}
where ${\mathcal{K}}$ is a tensor of the form
\begin{equation} \label{eqn-2D-tensor}
{\mathcal K} = \left[
\begin{array}{cc}
a & c \\
c & b
\end{array}
\right] ,~~~
\left\{
\begin{array}{l}
a = \cos^2(\theta) + \epsilon \sin^2(\theta) \\
b = \epsilon \cos^2(\theta) + \sin^2(\theta) \\
c = (1-\epsilon) \cos(\theta) \sin(\theta)
\end{array}
\right. .
\end{equation}
The parameter $0 < \epsilon \leq 1$ specifies the strength of anisotropy in the
problem, while parameter $\theta$ specifies the direction of anistropy. 
We use 2D domain $\Omega = [0,1] \times [0,1]$ with Dirichlet boundary conditions
with the following choice of problem parameters
\begin{equation*}
\begin{array}{ll}
\epsilon=1,~ d=10^4                &~~~ (x,y) \in {[0,1/2]} \times {[0,1/2]}, \\
\epsilon=1,~ d=0                   &~~~ (x,y) \in {(1/2,1]} \times {[0,1/2]}, \\
\epsilon=0.01,~ \theta=0,~ d=0     &~~~ (x,y) \in {[0,1/2]} \times {(1/2,1]}, \\
\epsilon=0.01,~ \theta=\pi/2,~ d=0 &~~~ (x,y) \in {(1/2,1]} \times {(1/2,1]}.
\end{array}
\end{equation*}
The discrete system is obtained by a finite element discretization of this problem with bilinear basis functions on uniform quadrilateral mesh.  The results of the algorithm for this test problem are provided in Figure~\ref{fig-cr-4reg};
we refer to~\cite{JBrannick_RFalgout} for additional 
numerical experiments and results obtained using CR in a classical AMG 
setup algorithm.  
In the reported results, an $F$-point Gauss Seidel 
variant of CR is used, the initial set $\coarsevar_0=\emptyset$ and the updating of $\coarsevar$ in the subsequent stages of the algorithm is based on adding an independent set of coarse-point candidates that are slow to converge in the CR process.
The precise implementation of the algorithm is given in~\cite{JBrannick_RFalgout}, see Algorithm 3.1 
in Section 3.  

Figure~\ref{fig-cr-4reg} illustrates the
coarsening stages and the final coarse grid chosen by the algorithm.  Because $\rho_{cr} 
> 0.7$, the algorithm enters the first coarsening stage.  Figure~\ref{fig-cr-4reg}(a) contains a plot of the
values of the candidate measures, scaled pointwise error values, that are used in selecting
coarse-grid candidates.  The candidate measures are large
in three regions of the domain, signaling that coarsening is needed there.  The
algorithm also correctly detects that coarsening is not needed in the bottom
left region of the domain where pointwise Gauss-Seidel is effective. The first
$C$-point set chosen is depicted by dotted squares in
Figure~\ref{fig-cr-4reg}(b).  For these updated sets, we
obtain $\rho_{cr} = 0.91$, which leads to a second stage.

The candidate measures at stage two of the algorithm are plotted in
Figure~\ref{fig-cr-4reg}(b) and the
algorithm signals that still further coarsening is needed only in the top two regions
of the domain.  The
result of the second stage is the coarse grid depicted in
Figure~\ref{fig-cr-4reg}(c).  The new CR convergence factor is $\rho_{cr}  = 0.49$,
so there are no additional stages in the algorithm.  Notice that the algorithm gives
semi-coarsening in $x$ in the top-left region where the problem is anisotropic
in $x$, semi-coarsening in $y$ in the top-right where the problem is anisotropic
in $y$, full-coarsening in the bottom-right where the problem is isotropic, and
no coarsening in the bottom-left where the mass term dominates.  Further, 
we mention that the algorithm does not make explicit use of an strength of connectivity 
meausre.

\begin{figure}[ht!]
\begin{center}
\begin{tabular}{@{}c@{}}
\includegraphics[width= 0.48\linewidth]{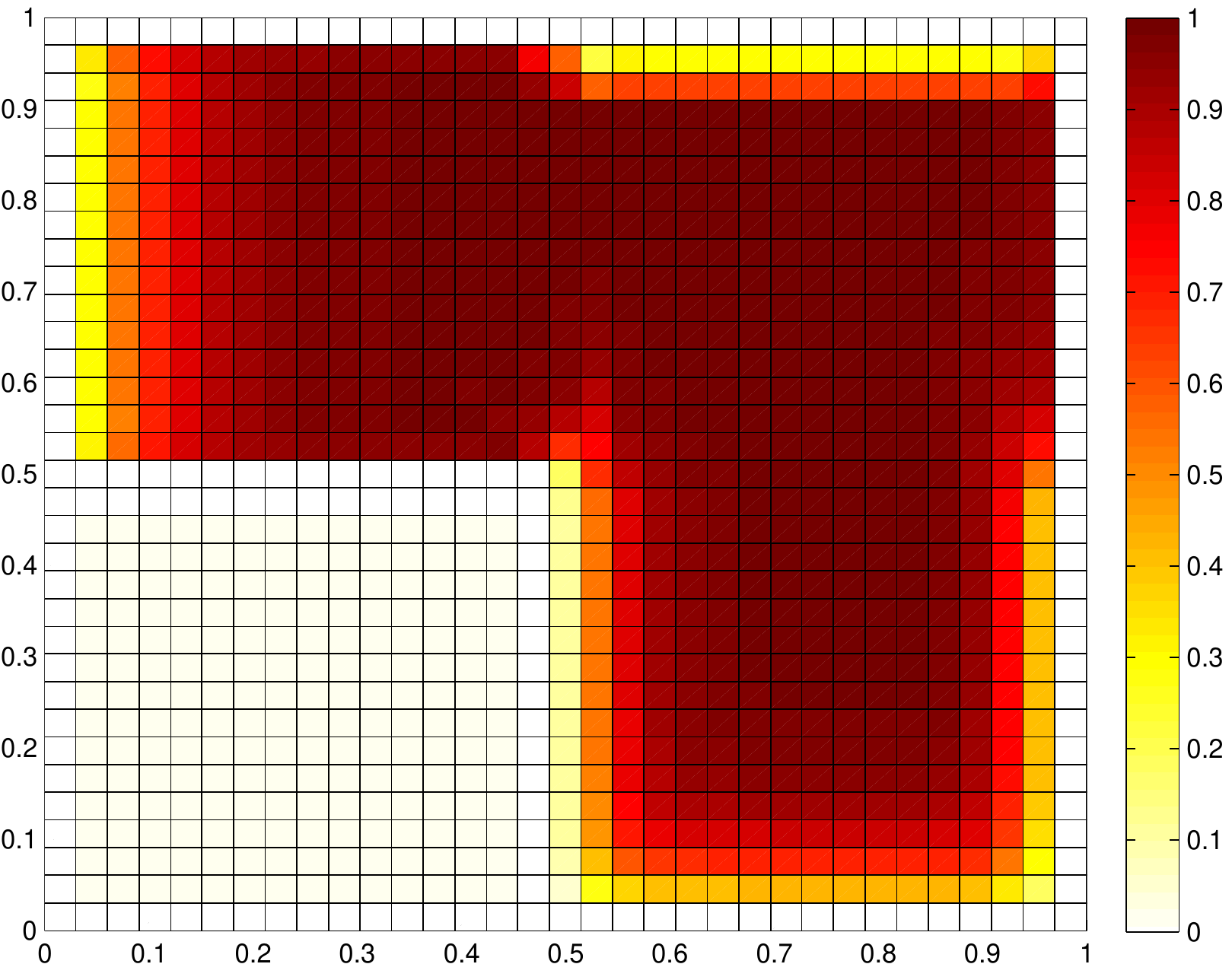} \\ (a)
\end{tabular}
\hfill
\begin{tabular}{@{}c@{}}
\includegraphics[width= 0.48\linewidth]{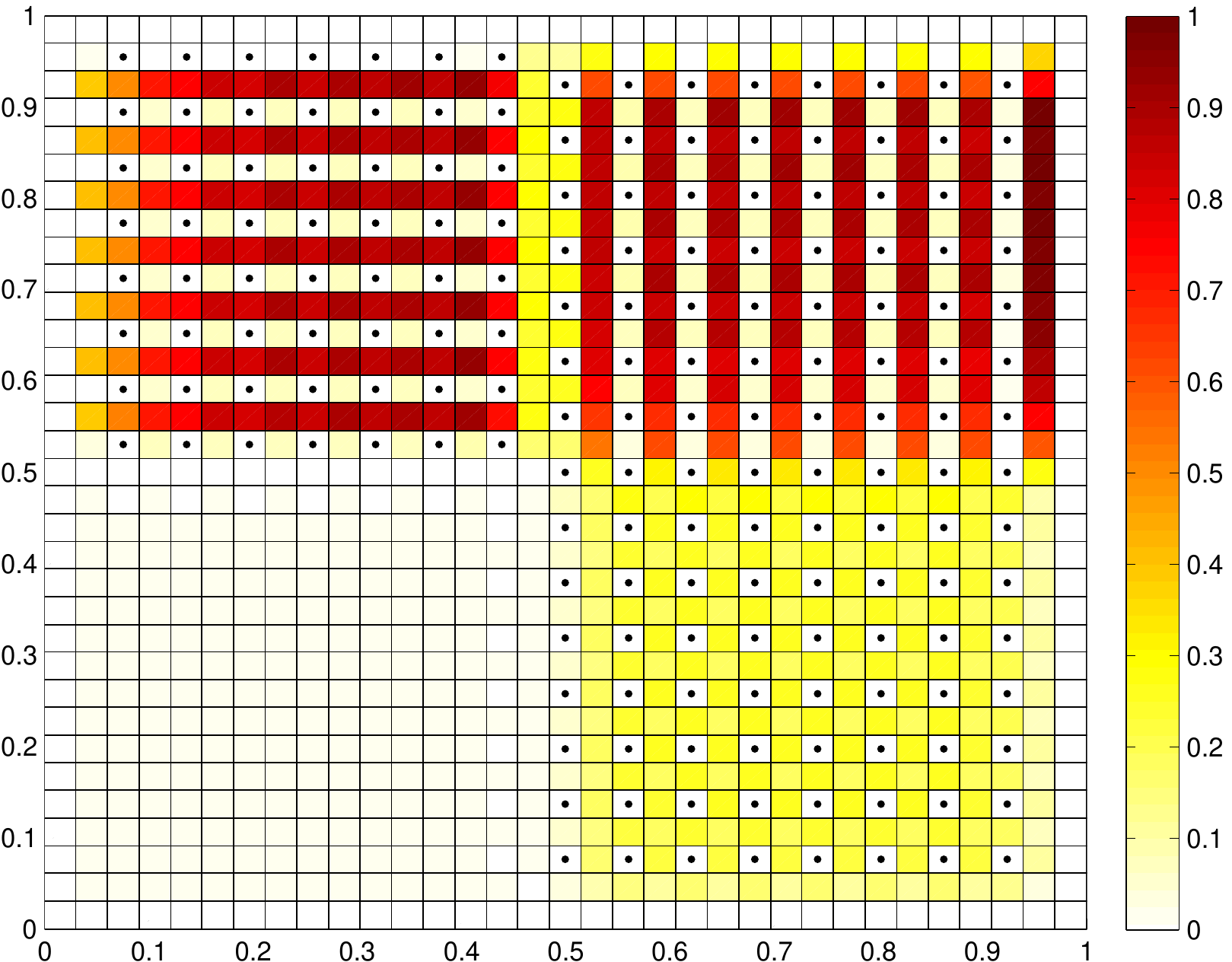} \\ (b)
\end{tabular}
\vspace{1em} \\
\begin{tabular}{@{}c@{}}
\includegraphics[width= 0.48\linewidth]{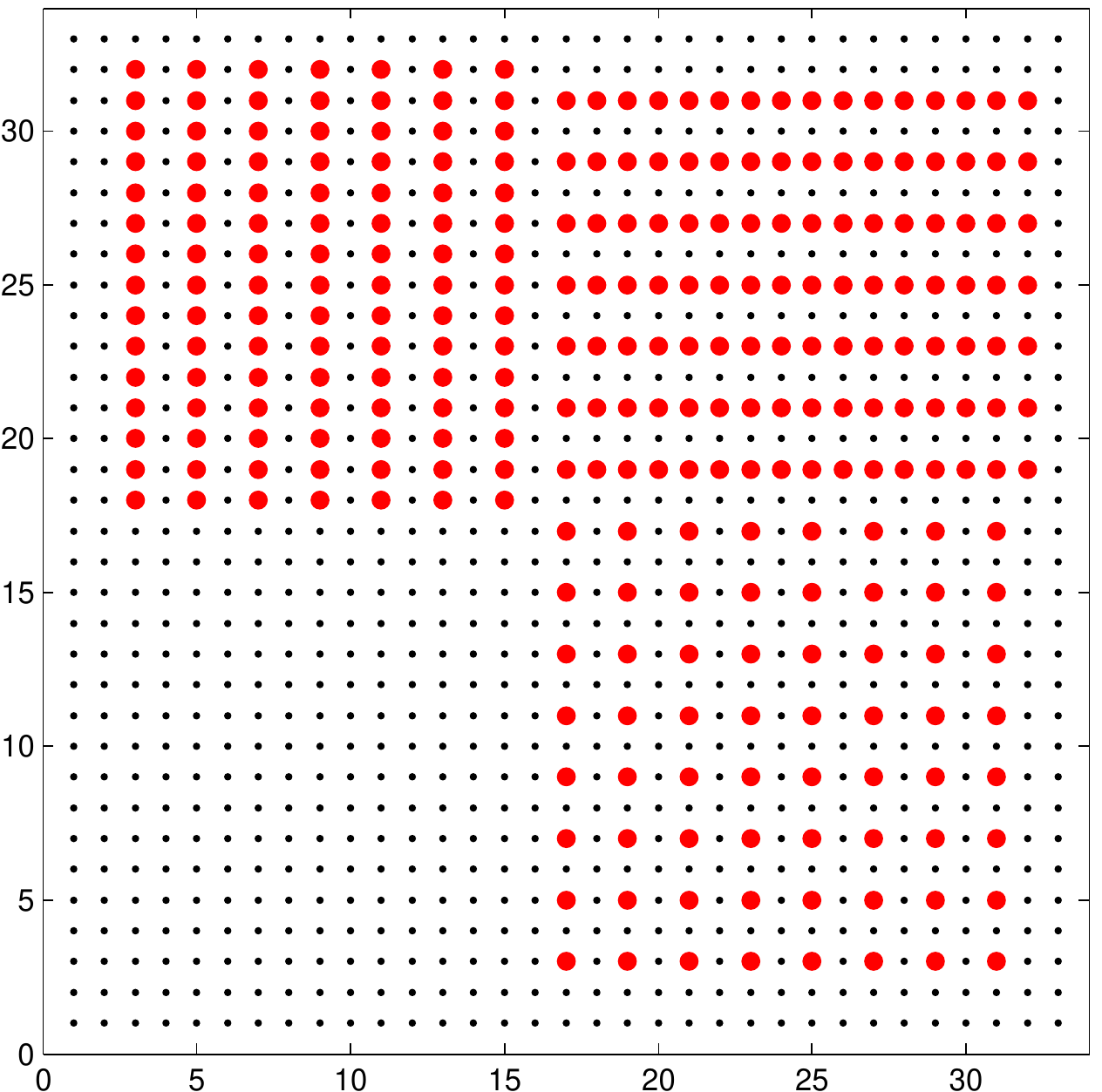}
\end{tabular}
\end{center}
\caption{
Candidate measures computed during
stage one (a) and two (b) for the 2D anisotropic diffusion problem.  Each square depicts the value
of the scaled pointwise error at grid point $i$ according to the color bar.  Dotted squares
indicate points currently in $\coarsevar$.  The final coarse grid is
given by the larger red circles in (c).}
\label{fig-cr-4reg}
\end{figure}
}

Generally, while CR algorithms have been successfully applied and developed for various model problems, there remain several design issues that must be addressed and remain open in developing a robust implementation
of the above CR algorithm.  
One issue that arises is formulating an appropriate relaxation scheme for the problem 
being solved.  If specific features of the problem are known, 
e.g., it comes from a system of PDEs, then an appropriate smoother
can be designed using knowledge of the problem.  More generally, 
the smoother can be adapted within the CR coarsening process.  For example,
in cases where $|\coarsevar|$ or $\rho_{cr}$ are too large, a block relaxation scheme can be used and the 
blocks can be adapted using CR to test the efficiency of the smoother.  
This is research topic which has not been explored much in the AMG setting, 
in particular, together with CR and algebraic distance.   {We mention the recent work
on adaptive smoothers in~\cite{Adap_Smooth} that shows the potential of such 
an approach.}

The two main research questions in CR that require further investigation and with which we conclude 
our discussion of CR are:

\begin{itemize}
\item  How to define an appropriate form of CR and related convergence rate estimate from $\rho_{cr}$ 
that accurately predicts the performance of the resulting two-grid method for a given interpolation strategy?
\item  How to choose $C_0$ and how to identify suitable local modifications to a given set of coarse variables
whenever the CR convergence rate is slow? 
\end{itemize}

The fact that fast CR rates imply the suitability of a given coarsening scheme has been established
theoretically~\cite{PanayotRob_2003} and by the extensive numerical experiments found in~\cite{JBrannick_RFalgout,Jacob}.  However, a precise relation between the convergence rate
of a given CR iteration and the resulting two-grid convergence is not known in general.  
Specific versions of CR have been proposed that do provide 
accurate predictions of the convergence rate of the associated two grid method, {\color{black}for example,
$F$-relaxation and} Habituated CR.  A single HCR iteration is given by relaxation steps applied 
to the entire fine-grid system, followed 
by imposing constraints on the resulting error restricted to the coarse variables:
$$e_C := Xe = 0,$$
where $X$ is a matrix that defines
the averaging process for deriving coarse variables.  
The solution to this problem can be calculated using the Kaczmarz iteration, starting with the error $e$ that results from the full relaxation steps as an initial guess or it
can be computed using the $\ell_2$ projection of $e$ on $\text{range} (X)$.  The latter approach
gives the overall HCR error propagation operator
\begin{equation}\label{eq:HCR}
 E = (I-\pi)(I - MA)^{\nu}, \quad \pi = X^H(XX^H)^{-1}X.
 \end{equation}
Unlike the $\finevar$-relaxation version of CR which can only be applied with the special coarsening used in classical AMG, and then only for simple relaxation schemes, HCR is applicable in general setting since it only involves
the fine-grid relaxation scheme followed by a simple projection based on the averaging matrix $X$.   Hence, in addition to giving a better qualitative prediction of the performance of the two-grid solver, HCR is more flexible and easier to process.

{\color{black}The sharpness of HCR in predicting the convergence rate
  of a two-grid method depends critically on the averaging matrix $X$
  used in enforcing the constraints and, as
shown in~\cite{JBrannick_RFalgout}, a suitable choice can be defined
by $X=P^H$, where $P$ denotes a sufficiently accurate and stable
interpolation operator. We note that this $P$ does not need to be the
same interpolation operator that is used in the AMG solver.
Generally, the accuracy of the iteration in~\eqref{eq:HCR} in
predicting the convergence rate of the two-grid method depends on the
energy norm of $\pi$, i.e., $\|\pi\|_A\geq 1$
(see~\cite{JBrannick_RFalgout}).  We mention, in addition, the work
in~\cite{SMacLYS07}, where it is shown that if the $C/F$ partitioning
is such that the weighted Jacobi $F$-relaxation converges quickly,
then a bound on the two-grid convergence rate follows.}

The issue of identifying suitable local modifications to the coarse variable set
whenever the CR convergence rate is too slow is a more subtle issue.  Assume 
that the relaxation scheme is fixed, $C$ is nonempty, 
and the convergence rate of CR is slow, indicating that additional variables need to be added 
to $C$.  
The usual approach is to find a candidate set, $\bar{F}$, of
$F$-points from which an independent set is added to the existing set
$C$.  
The primary issue is to determine a set of candidates $\bar{F}$.
The basic observation from~\cite{etna2000} that provides insight on this topic is that
a set of variables that should be added to the coarse variable set is exposed by the CR process; it is a subset of the set of variables which are slow to converge to zero when CR is applied. 
This observation motivated the CR algorithms developed and analyzed in~\cite{HCR,JBrannick_RFalgout}.

In~\cite{JBrannick_RFalgout}, a candidate set measure based on the values of the normalized pointwise
errors that result from the CR process was studied.  The set $\bar F$ is formed by adding to it $i\notin \coarsevar$ for which this error is large.   Intuitively, the approach is based on choosing additional coarse grid points in regions where the CR process is inefficient.  
In~\cite{HCR}, a similar algorithm was developed for the matrix $\tilde A$ that is first normalized such that each of its rows has $\ell_1$-norm equal to $1$.  
While the effectiveness of these strategies has been demonstrated by the numerous numerical experiments reported in~\cite{HCR,JBrannick_RFalgout},
there is still no rigorous understanding of or justification for using pointwise errors of CR to determine the variables that will most effectively reduce the convergence rate of the CR process when added to $C$.

Very recent work~\cite{BAMG_Aniso} on this topic has focused on incorporating the use of algebraic distance as a measure of strength of connectivity into the CR coarsening process.
The main idea in this context is to use algebraic distance to form a
strength of connectivity graph that identifies the coupling among
variables which can then be used in an independent set algorithm to
update $C$.
In the simplest form, the notion of algebraic distances is
{\color{black}straightforward: For any given pair of fine-grid variables $i,j \in \Omega$ compute
\begin{equation}\label{eq:AD}
 \mu_{ij}  =  \frac{
\sum_{\kappa = 1}^k\left(v_{i}^{(\kappa)} -  p_{ij} v_{j}^{(\kappa)}\right)^{2}} 
{\sum_{\kappa = 1}^k \bigg(v^{(\kappa)}_i\bigg)^2}, 
\end{equation}
where $p_{ij}$ is the coefficient that minimizes this least squares
problem for caliber-one interpolation of the form $v_i = p_{ij} v_j$.
Here, the set $\tvV = \{v^{(1)}, \ldots,
v^{(k)}\}$ consists of properly scaled smooth TVs, and the choice of
the denominator ensures that the expression is symmetric with respect
to the indices $i,j$. 
Statistically, the measure $\mu_{ij}$ is based on the idea that  
the nodes $i$ and $j$ in a graph are strongly connected only if $v_i$ and $v_j$ are highly correlated for all TVs. 

Another, more direct, way to write the algebraic distance measure is given by solving the LS problem \eqref{eq:AD} for the coefficient
$p_{ij}$ and substituting the minimizer into the definition of $\mu_{ij}$. 
The resulting 
$p_{ij} = \tfrac{\langle V_{(i)},V_{(j)}\rangle}{\langle V_{(j)},V_{(j)}\rangle}$ gives a rise to  an equation for the minimum of \eqref{eq:AD}:
\begin{equation}\label{eq:AD2}
  \mu_{ij} = 1- \frac{\big|\big\langle V_{(i)},V_{(j)}\big\rangle\big|^2}{\|V_{(i)}\|^2\|V_{(j)}\|^2}, 
  \end{equation}
  where 
with  $V = [v^{(1)} \: ... \:v^{(k)}] \in \mathbb{C}^{n\times k}$,
where    $\langle V_{(i)},V_{(j)}\rangle := \sum_{\kappa = 1}^k
v_i^{(\kappa)}v_j^{(\kappa)}$ measures the angle between $ V_{(i)}$
and $ V_{(j)}$.  From this explicit form of the minimum
  it is now easy to see that the measure $\mu_{ij}$ is symmetric with respect to $i,j$, as noted above.
}
{The algebraic distance measure is a specific case of the least squares interpolation process described in more detail in the next section.} 
More generally, the algebraic distance measure can be computed for sets of neighboring coarse points $\coarsevar_i$, using the more general 
LS problem defined in~\eqref{eq:LSfuncrowi}.  In this way, $\mu_{ij}$ can be used as an a posteriori measure of the quality of the interpolatory set $\coarsevar_i$ and the resulting row of interpolation.  
For cases where an algebraic notion of dependence is not sufficiently well defined, e.g., in the early stages of the BAMG setup process when the TVs may not be smooth enough to yield an accurate algebraic distance measure, one can assume any non-vanishing dependence to be strong in the CR coarsening process; this may result in slow coarsening, but then the 
scheme can be applied again in an iterative manner, to this new set instead of $C_0$.  The latter CR approach that does not rely on strength of connectivity to choose $\coarsevar$ was studied in the classical AMG setting in~\cite{JBrannick_RFalgout} and its extension that incorporates the algebraic distance measure in \eqref{eq:AD} was developed for anisotropic problems in~\cite{BAMG_Aniso}. 

Generally, such algorithms that use iterative strategies for choosing $C$ together with CR often lead to irregular alignment of coarse variables, even for structured grid problems, which in turn can result in an increase in the number of nonzeros in the coarse-grid matrix $A_c$.   This issue can be treated in latter stages of the BAMG process by using more aggressive coarsening, once the TVs become more smooth and the algebraic distance measure becomes more precise.  This is the case, for example, for strongly anisotropic problems where neighbors in the characteristic direction are much closer than in the other directions.  For such problems the appropriate long range connections in these characteristic directions are revealed only by very low energy TVs, which normally need to be the product of subsequent bootstrap setup cycles.
An algorithm that  both removes and adds variables to $C$ in each CR coarsening stage can be used to further control complexities, as needed.  We note that in our experience, for many problems
only the first coarsening stage in which $C_0$ is formed is needed, given sufficient accuracy of the algebraic distance measure.  

\section{Least squares interpolation}\label{sec:LS}
After the coarse grid $\coarsevar$, has been selected by compatible relaxation and algebraic distance, 
the interpolation operator is computed using a least
squares process. The least squares interpolation operator $P$ is defined to fit
collectively a set test vectors that should characterize the near kernel of the 
system matrix.   
Assuming the sets of interpolatory variables, $\coarsevar_i$, for each $i\in \finevar$, and a set of test vectors, $\tvV = \{v^{(1)}, \ldots,
v^{(k)}\}$, have been determined, the $i$th
row of $P$, denoted by $p_i$,
is defined as the minimizer of the local least squares problem: 
\begin{equation}\label{eq:LSfuncrowi}
\mathcal{L}(p_i) = 
\sum_{\kappa=1}^k\omega_k\left(v_{\{i\}}^{(\kappa)} - \sum_{j\in \coarsevar_{i}} \left(p_{i}\right)_{j} v_{\{j\}}^{(\kappa)}\right)^{2} \rightarrow \min.
\end{equation} 
Here, the notation $v_{\widetilde{\Omega}}$
denotes the canonical restriction of the vector $\widetilde{v}$ to
the set $\widetilde{\Omega} \subset \Omega$, e.g.,
$v_{\{i\}}$ is simply the $i$th entry of $v$.
Conditions on the uniqueness of the solution to minimization problem
\eqref{eq:LSfuncrowi} and an explicit form of the minimizer have been
derived in \cite{BAMG2010}.

The weights $\omega_{\kappa}>0$ are defined using the $A$-norms of the test vectors $\|v\|^2_A = \langle
Av,v\rangle$ when the system matrix $A$ is Hermitian and
positive definite (HPD)). 
 A main issue with this choice of weighting is that  
the LS process uses the entries of the test vectors restricted to 
{\em local} coarse neighborhoods to define the LS fit while the weights $\omega_{\kappa}$ are defined in terms 
of the {\em global} test vectors.  Hence, it may happen that a given test vector is globally smooth (i.e., 
produces a large weight and is emphasized in the LS process), but it does not provide an adequate local representation of smooth error.  This leads to a loss in accuracy of interpolation and slows down the convergence 
of the resulting two-grid method.  

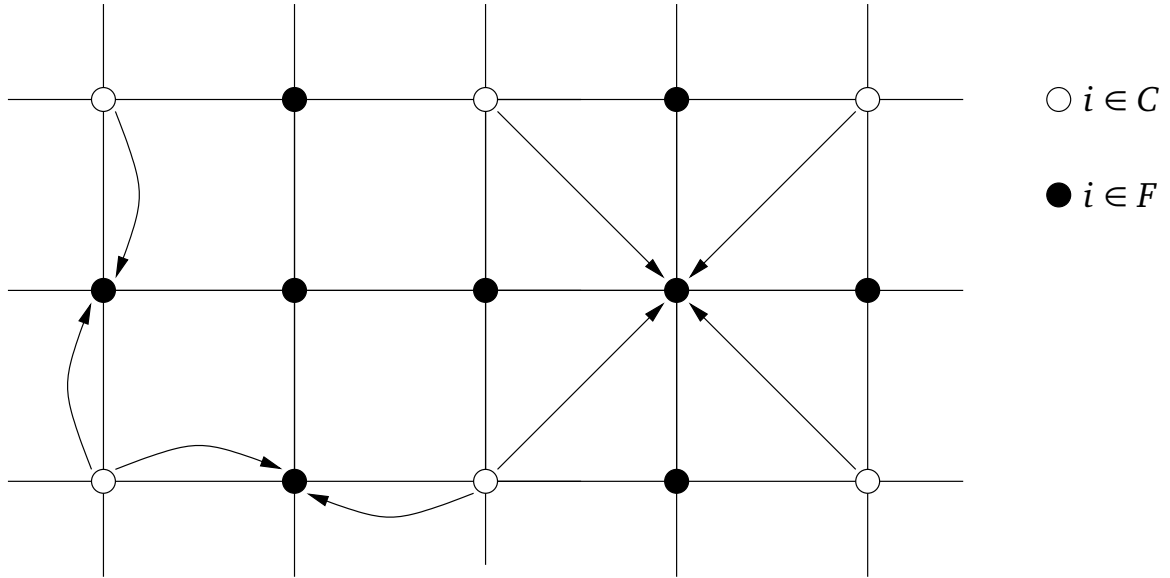
\begin{figure}
  \centering 
  \begin{picture}(0,0)%
\includegraphics{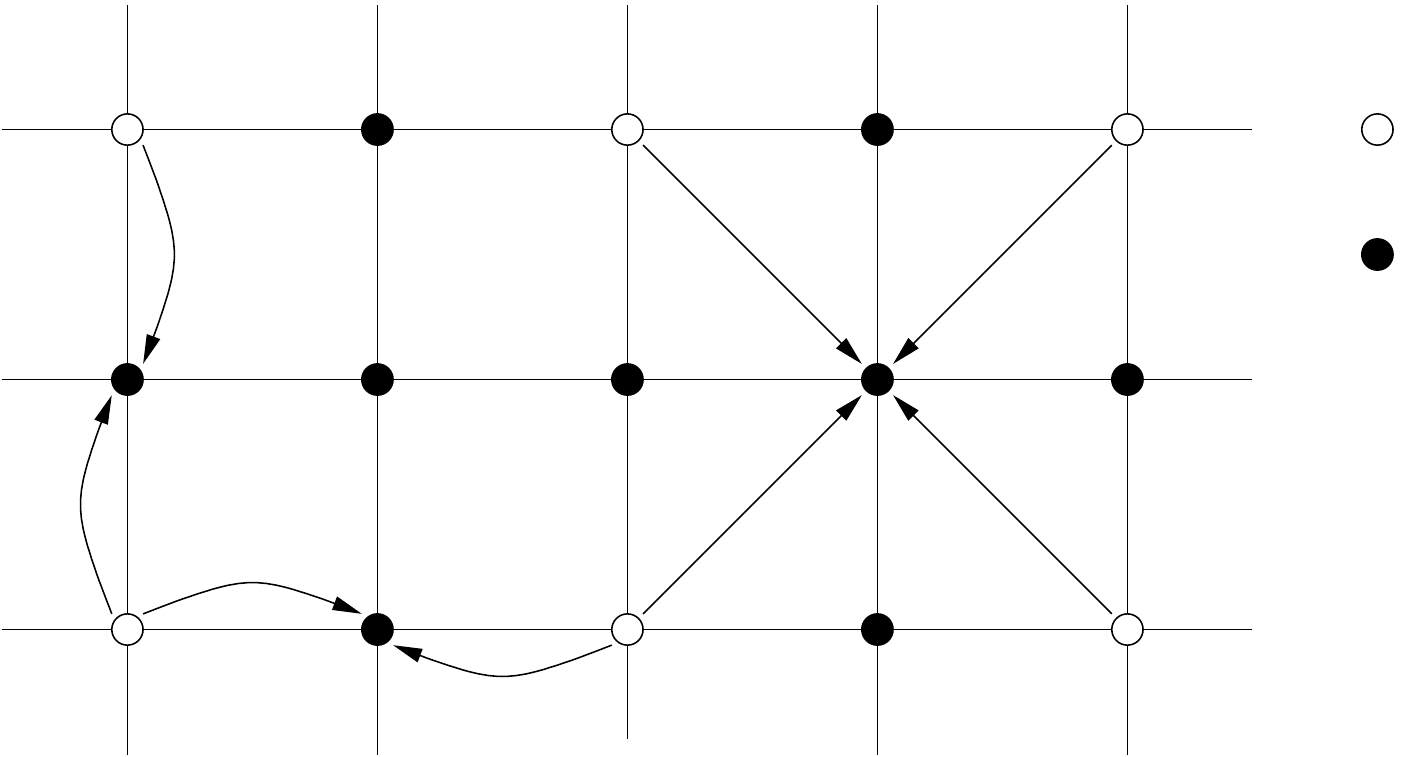}%
\end{picture}%
\setlength{\unitlength}{3947sp}%
\begingroup\makeatletter\ifx\SetFigFont\undefined%
\gdef\SetFigFont#1#2#3#4#5{%
  \reset@font\fontsize{#1}{#2pt}%
  \fontfamily{#3}\fontseries{#4}\fontshape{#5}%
  \selectfont}%
\fi\endgroup%
\begin{picture}(6777,3624)(2389,-3973)
\put(9151,-1636){\makebox(0,0)[lb]{\smash{{\SetFigFont{14}{16.8}{\rmdefault}{\mddefault}{\updefault}{\color[rgb]{0,0,0}$i \in \finevar$}%
}}}}
\put(9151,-1036){\makebox(0,0)[lb]{\smash{{\SetFigFont{14}{16.8}{\rmdefault}{\mddefault}{\updefault}{\color[rgb]{0,0,0}$i \in \coarsevar$}%
}}}}
\end{picture}%

\caption[Full coarsening and interpolation relations]{Full coarsening and interpolation relations for $i\in \finevar \setminus \coarsevar$.}
\label{fig:nr:scalar:interprel}
\end{figure}
{In~\cite{iBAMG} (see also~\cite{Brezina:2012:RCB})}, a residual-based LS form of BAMG interpolation (LSR) was introduced as a technique for 
addressing this issue of using global weights.  The basic idea in the LSR process
is to implicitly smooth the test vectors to produce an improved local representation 
of smooth error.  In this respect, the approach is related to the adaptive relaxation scheme described by Brandt in~\cite{BAMG1} (and suggested in~\cite{etna2000}) that applies relaxation selectively to points exhibiting especially large residuals, as they occur for example in problems with singularities.
Specifically, in LSR interpolation the test vectors are (temporarily) updated as follows:
\begin{equation}\label{eq:LSfsmooth}  v_{\{i\}}^{(\kappa)} \leftarrow  v_{\{i\}}^{(\kappa)} - \frac{1}{d_{ii}}  r_{\{i\}}^{(\kappa)}, \end{equation}
and then used in computing $p_i$, again as the solution to \eqref{eq:LSfuncrowi}.  

Table~\ref{tab:nr:lsinterpketaN01} contains results from~\cite{BAMG2010} for a variational two-grid method resulting from using standard (full) coarsening and interpolation with $|\coarsevar_i| \leq 4$ in \eqref{eq:LSfuncrowi} (see Figure~\ref{fig:nr:scalar:interprel}) for the original LS formulation and the LSR version
applied to the central finite difference discretization of the Poisson problem with homogenous Dirichlet boundary conditions on a uniform quadrilateral grid.  

The reported results are the approximate asymptotic convergence rates for different values of the number relaxed test vectors, $k_r$, and the number of Gauss Seidel iterations, $\eta$, used to compute them.  {Note that here the total number of TVs is set to $k=k_r$, whereas, when a multilevel eigensolver is combined with relaxation
{\color{black}(an approach that we consider in the next section)}
the total number of TVs is given by $k=k_r+k_e$, with $k_e$ the number of TVs computed by the multilevel eigensolver described in the following section.}  The test vectors are computed by applying Gauss Seidel iterations with lexicographic ordering to the Laplace problem starting with $k$ distinct initial guesses, generated randomly with a normal distribution with expectation zero and variance one, $N(0,1)$, and the update in \eqref{eq:LSfsmooth} used in the LSR formulation is applied only to $20\%$ of the entries of the TVs for which the associated values of the residual $r_{\{i\}}^{(\kappa)}$ are largest in absolute value.

\begin{table}
  \begin{center}
    \begin{tabular}{|c||c|c|c|c|c|c|c|}
      \hline
     $k$ & $6$ & $7$ & $8$ & $10$ & $12$\\ 
             \hline 
        \hline 
       $\eta =2$ & $.934$    $(.866)$  &  $.916$    $(.827)$  &  $.904$ $(.806)$ &   $.878$    $(.786)$    & $.861$    $(.759)$  \\
      $\eta =4$ & $.735$    $(.524)$ &   $.681$    $(.453)$  &  $.648$    $(.403)$  &  $.608$    $(.269)$ & .575 (.316) \\
      $\eta =6$ & $.522$    $(.294)$   & $.451$    $(.242)$  &  $.425$    $(.205)$   & $.381$    $(.154)$ &.305 (.136) \\
      $\eta =8$ & $.372$    $(.169)$  &  $.308$    $(.140)$  &  $.287$    $(.123)$   & $.238$   $(.054)$  & .216 (.061)\\
      \hline 
    \end{tabular}
    \caption[Two-grid LS interpolation forPoisson probelm ($k$ vs.\ $\eta$)]{Asymptotic convergence rate estimates, $\rho$, of the two-grid methods with two pre- and two post-smoothing steps constructed using LS (residual-based LS) interpolation for the standard finite difference discretization of the Poisson problem with homogenous Dirichlet boundary conditions on a uniform rectangular grid with $h=1/64$.  
     }\label{tab:nr:lsinterpketaN01} 
  \end{center}
\end{table} 

The improved performance of the method that results from increasing $\eta$
is perhaps expected, since applying more relaxation sweeps to the test vectors should
give a better approximation of smooth error.  The situation gets more interesting and more  complicated when we fix the number of Gauss Seidel iterations and then increase the sampling space (i.e, vary the number of test vectors), which results in an ad hoc (randomized) approach for sampling the dominant part of the space that is not significantly reduced by the smoother. 
Recently, randomized methods for range approximation
have been given particular attention in data mining applications and further studies 
of the LS process along these lines are still needed.  {The paper~\cite{antrg} summarizes the latest
results obtained on this topic and highlights the gains which are possible with judicious modifications of the
ad hoc approach. }

\begin{table}[!ht]
  \begin{center}
    \begin{tabular}{|c|c|c|c|c|c|c|}
      \hline
      $h$ & $1/32$ & $1/64$ & $1/128$ & $1/256$ & $1/512$ \\
      \hline 
    $\rho$ & $.267$ $(.075)$ & $.648$ $(.403)$ & $.886$
      $(.758)$ & $.966$ $(.917)$ & $.990$ $(.977)$\\
            \hline
    \end{tabular}
    \caption[Two-grid LS interpolation ($k_{r}=8, \eta=4$) for FD Laplace]{Asymptotic convergence rate estimates, $\rho$, of the LS (LSR) based two-grid methods with two pre- and two post-smoothing steps applied to the same FD Laplace problem as in Table 1.  The solver is constructed using the LS (LSR) schemes with $\eta=4$ and $k_{r}=8$.\label{tab:nr:lsinterp84N01}}
     \end{center}
\end{table} 

In Table~\ref{tab:nr:lsinterp84N01}, we report results from~\cite{BAMG2010} for a two-grid method applied to the same central finite difference Poisson problem for different choices of the grid spacing $h$.  The solver is constructed using the LS and LSR formulations for defining $P$ with a fixed number of relaxed test vectors $k_{r} = 8$ and smoothing iterations
$\eta = 4$ independent of the choice of $h$.  Here, we see that the convergence rates of both the LS and LSR based two-grid methods increase
as the size of the problem is increased, suggesting that applying a constant number of relaxation steps to compute a fixed number of test vectors does not produce a sufficiently accurate local representation of  the algebraically smooth error for defining least squares interpolation.  {The numerical results reported in~\cite{iBAMG,Brezina:2012:RCB}
further support this claim.}

In Table~\ref{tab:nr:lsinterp714N01}, we present results of the two-grid solvers obtained by adding the constant vector, ${\bf 1}$, to $\mathcal{V}$ and again applying the LS and LSR schemes to the FD Poisson problem.  Here a fixed number of
test vectors $k_{r} = 7+1$ are used and $\eta =4$ relaxation steps are used to compute them.  
\begin{table}[!ht]
  \begin{center}
    \begin{tabular}{|c|c|c|c|c|c|c|}
      \hline
      $h$ & $1/32$ & $1/64$ & $1/128$ & $1/256$ & $1/512$ \\
      \hline 
       $\rho$ & $.089$ ($.039$) & $.121$ ($.040$) & $.130$
      ($.042$) & $.148$ ($.042$) & $.150$ ($.043$)\\
      \hline
    \end{tabular}
    \caption[Two-grid LS interpolation ($k_{r}=7+\mathbf{1}, \eta=4$) for FD Laplace]{Asymptotic convergence rate estimates, $\rho$, of the LS (LSR) based two-grid methods with two pre- and two post-smoothing steps applied to the same FD Laplace problem as in Table 1.  The solver is constructed using the LS (LSR) schemes using $\eta=4$ Gauss Seidel iterations to compute the seven initially random test vectors; the additional test vector is the constant vector chosen a priori.
      \label{tab:nr:lsinterp714N01}}
    \end{center}
\end{table} With this modification the performance of the LS and LSR based
two-grid solvers improve significantly.  In fact, both the LS and LSR methods perform nearly optimally for the problem sizes
considered.  This clearly demonstrates the ability of the LS process to incorporate geometric and spectral information from the problem when it is available.  We note, in addition, that with our choice of full coarsening and interpolation relations (see Figure~\ref{fig:nr:scalar:interprel}), the
resulting interpolation is not derived from immediate neighbors, which 
further demonstrates the flexibility of the least squares formulation in that it naturally 
accommodates long-range interpolation.
We refer to the paper~\cite{BAMG_Aniso} in which algebraic distance and long-range LS interpolation were studied for
non-grid-aligned anisotropic problems for additional details and results on this topic.  

The LSR process provides a simple technique for improving the accuracy of interpolation in cases where the global weights are large but the TVs are not sufficiently smooth in some local regions.  As the numerical experiments show, the main point is to obtain extra smoothness (i.e., smaller residuals) at the points
where interpolation is defined.  And such large residuals in general are the result of special local configurations in the random vectors used as initial guesses.  This issue can also be inexpensively treated using adaptive relaxation together with several short $V$-cycles based on the 
interpolation that is constructed for the first few levels.  

The related problem of determining appropriate local weights in the LS definition of interpolation is more 
complicated and has not been resolved in a general algebraic setting, except for finite element discretizations of elliptic
boundary value problems.  For FE systems, the local weights can be computed using the local stiffness matrices.  
The problem of determining local weights is also related to the question
of deriving local equations for computing the test vectors, as in spectral AMG~\cite{rhoAMGe}.

\section{A bootstrap algorithm for computing the test vectors}\label{sec:BS}
\begin{figure}
\begin{center}
     \tikzstyle{greenpoint}=[circle,inner sep=0pt,minimum size=2mm,draw=black!100,fill=black!30]
    \tikzstyle{whitepoint}=[rectangle,inner sep=0pt,minimum size=2mm,draw=black!100,fill=black!80]
    \tikzstyle{redpoint}=[diamond,inner sep=0pt,minimum size=2.55mm,draw=black!100,fill=black!80]
    \tikzstyle{blackpoint}=[circle,inner sep=0pt,minimum size=2mm,draw=black!100,fill=black!100]
    \tikzstyle{bluepoint}=[circle,inner sep=0pt,minimum size=2mm,draw=black!100,fill=black!0]
    \resizebox{\textwidth}{!}{\begin{tikzpicture}
      \draw [sharp corners] (0,0) node[blackpoint] {} --
      ++(300:1cm) node[blackpoint] {} --
      ++(300:1cm) node[blackpoint] {} --
      ++(300:1cm) node[blackpoint] {} --
      ++(300:1cm) node[redpoint] (L4) {} --
      ++(60:1cm) node[bluepoint] (L3) {} --
      ++(60:1cm) node[bluepoint] (L2) {} --
      ++(60:1cm) node[bluepoint] (L1) {} --
      ++(60:1cm) node[whitepoint] (L0) {} --
      ++(0:1cm) node[greenpoint] {} --
      ++(300:1cm) node[greenpoint] (end) {};
      \draw[dashed] (end) --  ++(300:1cm);
      
      \draw (L0) + (3cm,0) node[blackpoint,label=right:{\small Relax on $Av=0, v \in \mathcal{V}^r$, compute $P$}] {};
      \draw (L1) + (3.5cm,0) node[redpoint,label=right:{\small Compute $v$, s.t., $Av=\lambda Tv$, update $\mathcal{V}^e$}] {};
      \draw (L2) + (4.0cm,0) node[bluepoint,label=right:{\small Relax on $Av = \lambda T v, v \in \mathcal{V}^e$}] {};
      \draw (L3) + (4.5cm,0) node[greenpoint,label=right:{\small Relax on and / or solve using MG $Av=0, v \in \mathcal{V}^r$ and $Av = \lambda T v, v \in \mathcal{V}^e$, recompute $P$}] {};
      \draw (L4) + (5cm,0) node[whitepoint,label=right:{\small Test
        MG method, update $\mathcal{V}$}] {};

      \draw (L3) +(0,4cm) node (start) {};
      \draw[sharp corners] (start) ++(300:-3cm) node[blackpoint] {} --
      ++(300:1cm) node[blackpoint] {} --
      ++(300:1cm) node[blackpoint] {} --
      ++(300:1cm) node[redpoint] {} --
      ++(60:1cm) node[greenpoint] {} --
      ++(300:1cm) node[redpoint] {} --
      ++(60:1cm) node[bluepoint] {} --
      ++(60:1cm) node[greenpoint] {} --
      ++(300:1cm) node[greenpoint] {} --
      ++(300:1cm) node[redpoint] {} --
      ++(60:1cm) node[greenpoint] {} --
      ++(300:1cm) node[redpoint] {} --
      ++(60:1cm) node[bluepoint] {} --
      ++(60:1cm) node[greenpoint] {} --
      ++(300:1cm) node[greenpoint] {} --
      ++(300:1cm) node[redpoint] {} --
      ++(60:1cm) node[bluepoint] {} --
      ++(60:1cm) node[bluepoint] {} --
      ++(60:1cm) node[whitepoint] {} -- 
      ++(0:1cm) node[greenpoint] {} --     
      ++(300:1cm) node[greenpoint] {} --
      ++(300:1cm) node[greenpoint] {} --
      ++(300:1cm) node[redpoint] {} --
      ++(60:1cm) node[greenpoint] (end) {};
      \draw[dashed] (end) --  ++(300:1cm);
    \end{tikzpicture}}
  \caption{Galerkin Bootstrap AMG $W$-cycle and $V$-cycle setup schemes.\label{fig:boot:setupcycle}}
\end{center}
\end{figure}
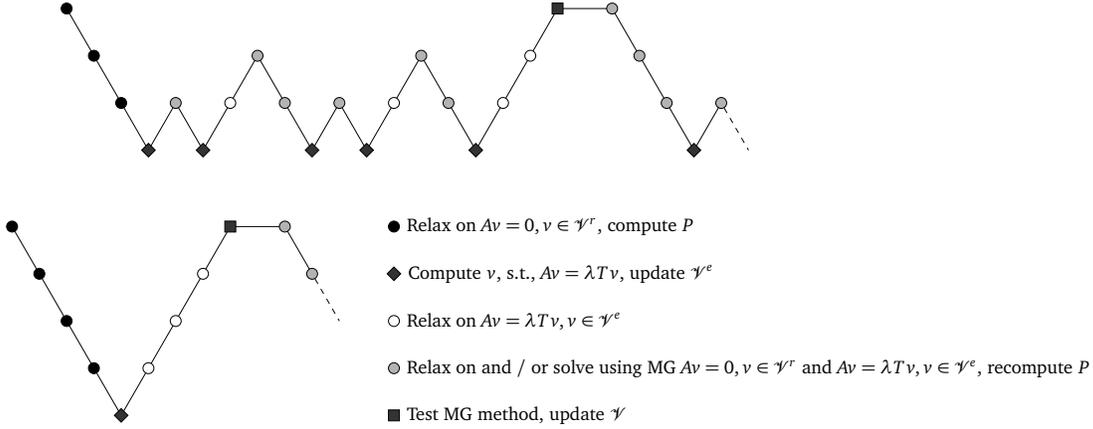

The BAMG process provides a general approach for computing the TVs used in 
algebraic distance and LS interpolation.  For definiteness, we provide an overview 
of one specific version of bootstrap AMG designed and analyzed in~\cite{BAMG2010}. 
The algorithm begins with 
relaxation applied to the homogenous system, 
\begin{equation}\label{homog}
A_lx_l =0,
\end{equation}
on each grid, $l = 0,..., L-1$; assuming that a priori knowledge of the 
algebraically smooth error is not available, these vectors are initialized randomly on the finest grid, whereas on all coarser grids
they are defined by restricting the existing test vectors computed on the previous finer grid.  
Given interpolation constructed using CR, algebraic distances and the LS process, 
the coarse-grid operators are computed, e.g., using the variational definition when $A$ is HPD.  
Once an initial MG hierarchy has been computed, 
the current sets of TVs are further enhanced on all grids using the 
existing multigrid structure.  Specifically, the given hierarchy is used to formulate a multigrid eigensolver
which is then applied to an appropriately chosen generalized eigenproblem to compute additional test vectors.  
{We note that the use of a coarse-grid eigensolver to compute test vectors has been studied theoretically in the context
of reduction-based AMG in~\cite{MacRed} and in the context of adaptive smoothed aggregation
in~\cite{GESA}.  The approach outlined here differs in its use of LS interpolation together with 
a multilevel eigensolver as considered in~\cite{BAMG2010}.}
This overall process is then repeated with the current AMG method applied in addition to (or replacing)
relaxation as the solver for the homogenous systems in \eqref{homog}.
Figure \ref{fig:boot:setupcycle} provides an schematic outline of the bootstrap $V$- and $W$-cycle setup algorithms.  {In general, $V^m$ and $W^m$ denote setup algorithms that use $m$ iterations of the $V$- and $W$-cycles, respectively.}

The rationale behind the multilevel generalized eigensolver (MGE) is as follows.  Assume an initial multigrid hierarchy has been constructed.  Given the initial Galerkin operators $A_{0}, A_{1}, \ldots, A_{L}$ on each level
and the corresponding interpolation operators $P_{l+1}^{l}, l =
0,\ldots,L-1$, define the composite interpolation operators as
$P_{l} = P_{1}^{0}\cdot \ldots \cdot P_{l}^{l-1},\ l = 1, \ldots, L$.  Then, 
for any given vector $x_{l} \in \mathbb{C}^{n_l}$ we have
$
  \innerprod[A_{l}]{x_{l}}{x_{l}} = \innerprod[A]{P_{l}x_{l}}{P_{l}x_{l}}.
$  
Furthermore, defining $T_{l} = P_{l}^{H}P_{l}$ we obtain
\begin{equation*}
  \frac{\innerprod[A_{l}]{x_{l}}{x_{l}}}{\innerprod[T_{l}]{x_{l}}{x_{l}}} = \frac{\innerprod[A]{P_{l}x_{l}}{P_{l}x_{l}}}{\innerprod{P_{l}x_{l}}{P_{l}x_{l}}}.
\end{equation*} 
This observation in turn implies that 
on any level $l$, given a vector $x^{(l)} \in \mathbb{C}^{n_l}$ and $\lambda^{(l)}\in
  \mathbb{C}$ such that
$  A_l x^{(l)} = \lambda^{(l)}T_l x^{(l)},
$
we have Rayleigh quotient (RQ)
\begin{equation}\label{eq:evalapprox}
  \mbox{rq}(P_l x^{(l)}) :=
  \frac{\innerprod[A]{P_{l}x^{(l)}}{P_{l}x^{(l)}}}{\innerprod{P_{l}x^{(l)}}{P_{l}x^{(l)}}}
  = \lambda^{(l)}.
\end{equation}
This provides a relation among the eigenvectors and eigenvalues
of the operators in the multigrid hierarchy on all levels with 
the eigenvectors and eigenvalues of the finest-grid operator $A$
that can be used to measure the accuracy of interpolation for smooth
modes.  
Note that the eigenvalue approximations in~\eqref{eq:evalapprox} are continuously updated within the algorithm so that the overall approach resembles an inverse Rayleigh-Quotient
iteration found in eigenvalue computations (cf.~\cite{JWilkinson_1965}).  
The overall MGE process is motivated by the following observation.

Consider the generalized eigenvalue problem
  \begin{equation}\label{eq:eva}
    A_{j}x_{j} = {\color{black}\lambda^{(j)}}T_{j}x_{j},
  \end{equation} on two subsequent grids, $j=l, l-1$, with 
  $\ {\color{black}\lambda^{(l)}}$ denoting the RQ of a given approximation of an eigenvector on the coarser grid and $ {\color{black}\lambda^{(l-1)}}$ denoting the RQ of the vector obtained by applying relaxation to the system~\eqref{eq:eva} with this vector interpolated to the next finer grid as the initial guess.  Define the \emph{eigenvalue approximation measure} $\tau^{(l,l-1)}_{\lambda}$ by
    \begin{equation}\label{eq:EWapproxmeasure}
  {\color{black}   \tau^{(l,l-1)}_{\lambda} =
    \frac{|\lambda^{(l)}-\lambda^{(l-1)}|}{|\lambda^{(l-1)}|}}.
  \end{equation}
Then, at any stage of the MGE iteration, a large value of $\tau^{(l,l-1)}_{\lambda}$ indicates that the hierarchy should be recomputed to incorporate this relaxed eigenvector approximation.  Otherwise, relaxation has not significantly changed this vector and it is accurately represented by the existing interpolation operator $P_l^{l-1}$.  In this way, the MGE serves as an approach for efficiently identifying components that must  be interpolated accurately (e.g., the low modes of $A$), and with the use of the eigenvalue approximation measure a technique for determining if the current $P$ approximates them sufficiently well.  
For additional details of the algorithm and its implementation we refer to the paper~\cite{BAMG2010}.  

To illustrate the effect of the MGE in the bootstrap process, we provide results from~\cite{BAMG2010}, in which a $V^2$-cycle bootstrap setup using four pre- and post-smoothing steps is applied to compute the set of relaxed vectors $\mathcal{V}^r$ and set of bootstrap vectors $\mathcal{V}^e$ coming from the MGE process, with $k_r=|\mathcal{V}^r| = 8$ and $k_e= |\mathcal{V}^e|=8$, respectively.  The sets $\mathcal{V}^r$ and $\mathcal{V}^e$ are then combined to form the set of TVs $\mathcal{V}$ that is used to compute the least squares interpolation operator on each level.   
In these tests, only relaxation is applied to the homogenous systems in both setup cycles to update the sets $ \mathcal{V}^r$.
The coarse grids and sparsity structure of interpolation are defined as in the previous tests on all levels
and the problem is coarsened to a coarsest level with $h=1/16$.
  
 As the results reported in Table~\ref{tab:nr:lsinterpMG84N01V2} show, using the MGE to enhance the 
test vectors consistently improves the performance of the resulting solvers
when compared to the results reported in Table~\ref{tab:nr:lsinterp84N01},
especially for the LSR scheme.
\begin{table}[!ht]
  \begin{center}
    \begin{tabular}{|c|c|c|c|c|c|c|}
      \hline
      $h$ & $1/32$ & $1/64$ & $1/128$ & $1/256$ & $1/512$ \\
      \hline 
     $\rho$ & $.041$ ($.038$) & $.062$ ($.041$) & $.075$
      ($.043$) & $.125$ ($.043$) & $.971$ ($.043$)\\
      \hline
    \end{tabular}
    \caption[$V^{2}$-BootAMG ($k_{r}=8, \eta=4$) for FD Laplace]{
   Asymptotic convergence rate estimates, $\rho$, of the LS (LSR) based $V(2,2)$ BAMG solvers applied to the FD Laplace problem.  The BAMG solvers are constructed using a $V^2$-cycle setup with $\eta=4$ and $k_{r}=k_e=8$.  
\label{tab:nr:lsinterpMG84N01V2}}
  \end{center}
\end{table} 
Here, the LSR approach yields an optimal method, whereas for the LS scheme, as the problem size increases the convergence rate of the solver also increases substantially -- from $\rho \approx.125$ to $\rho \approx.971.$ 

Next, we consider results of the $W$-cycle setup, illustrated at the top of Figure~\ref{fig:boot:setupcycle}.
In Table~\ref{tab:nr:lsinterpMG84N01W}, we present results
for this scheme again with $k_{r}=k_e=8$ and using only relaxation to solve the homogenous systems in both setup cycles to update the sets $ \mathcal{V}^r$.  We notice a marked improvement in the performance of the LS-based solver constructed using a $W$-cycle setup over that of the solver
constructed using a $V$-cycle setup and note that the LSR scheme again scales optimally for the given problem sizes.  The improved performance of the $W$-cycle setup indicates the effectiveness of using additional  coarse level processing in computing accurate TVs. 

\begin{table}[ht!]
  \begin{center}
    \begin{tabular}{|c|c|c|c|c|c|c|}
      \hline
      $h$ & $1/32$ & $1/64$ & $1/128$ & $1/256$ & $1/512$ \\
      \hline 
    $\rho$ & $.042$ ($.038$) & $.062$ ($.041$) & $.073$
      ($.043$) & $.130$ ($.043$) & $.161$ ($.044$)\\
      \hline
    \end{tabular}
    \caption[$W$-BootAMG ($k=j=8, \eta=4$) for FD Laplace]{Asymptotic convergence rate estimates, $\rho$, of the LS (LSR) based $V(2,2)$ BAMG solvers applied to the FD Laplace problem.  The BAMG solvers are constructed using a $W$-cycle setup (see~Figure~\ref{fig:boot:setupcycle}) with $ \eta=4$ and $k_{r}=k_e =8$.\label{tab:nr:lsinterpMG84N01W}}
  \end{center}
\end{table}

\begin{figure}[!t]
\centering{
\includegraphics[width=.65\textwidth]{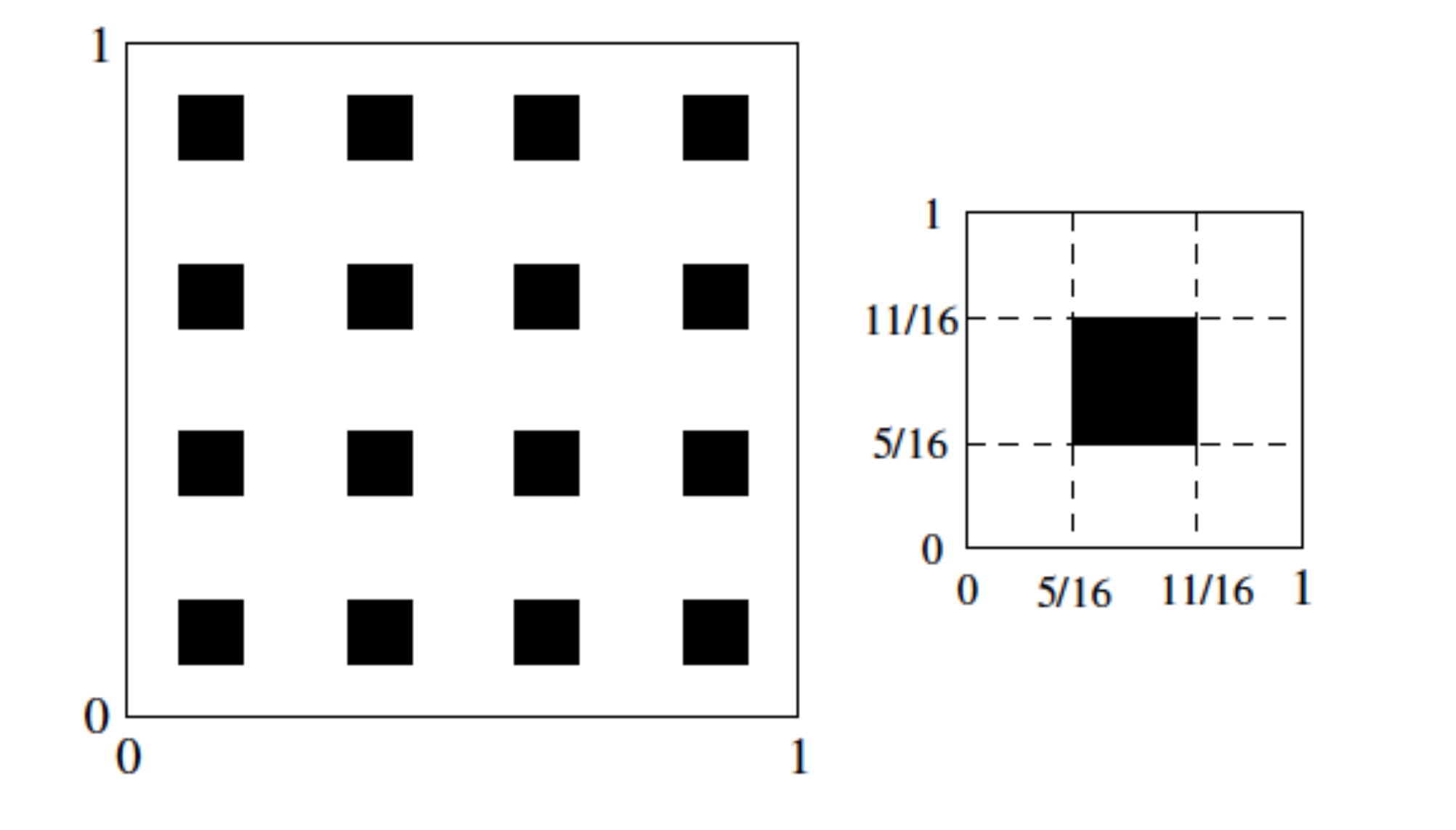}   }
\caption{Tiling of the domain for the interface problem \eqref{eq:jumps}.  The figure on the left corresponds to $\tau = 4$ and the one on the right to $\tau =1$. \label{fig:jumps}}  
\end{figure}

We conclude the numerical experiments with results of the BAMG setup algorithm applied to the 
same bilinear finite element discretization of the 2D diffusion problem considered earlier
\begin{equation}
 - \nabla \cdot \mathcal{K}(x,y) \nabla u(x,y) = f(x,y),
 \end{equation}
 with homogenous Dirichlet boundary conditions.  
However, now we assume that $\mathcal{K}(x,y)$ denotes the two-scale permeability:
$$\mathcal{K} = \bigg\{ \begin{array}{lr}
1 & (x,y) \in \Omega_1  \\
10^\nu  & (x,y) \in \Omega_2
\end{array},$$
with $\Omega_1$ given by the union of the black subdomains in Figure~\ref{fig:jumps}
and $\Omega_2 = \Omega \setminus \Omega_1$.  
The problem is discretized on a uniform quadrilateral grid with various choices of the grid spacing
$h$ so that the grid aligns with the interfaces
for the given problem sizes (see the plot on the right of Figure~\ref{fig:jumps}).  We consider the values $\nu=-2,-8$ 
and various tilings, denoted by $\tau=1,4,8,16$ (see Figure~\ref{fig:jumps}) in testing the BAMG algorithm.

The BAMG method uses standard (full) coarsening with nearest-neighbor interpolation on all levels of the hierarchy, and the problem is coarsened until the coarsest systems has grid spacing $h=1/16$.  We note that on coarse levels, the grids and the interpolation patterns do not align with the jumps at the interfaces of $\Omega_1$ and $\Omega_2$.
We consider tests of the $W^2$-cycle setup with $\eta =4$ relaxation steps to compute the $k_{r}=k_e=8$ TVs.   For the results reported in the bottom two tables of Table~\ref{tab:aBAMG} we use the solver resulting from the first $W$-cycle setup in an adaptive step in which we solve the homogenous system $Ax=0$ using five $V(2,2)$ cycles on the finest level only.  
Here, the initial guess is the TV from $ \mathcal{V}^r$ with the smallest energy, which is then replaced
by the resulting approximation as the new TV in $ \mathcal{V}^r$.

\begin{table}[!ht]
    \begin{center}\small
      \begin{tabular}{lll}
        \begin{tabular}{|c|c|c|c|c|c|c|c|}
          \hline
          $\tau$\\
          \hline
          $1$\\
          \hline
          $4$\\
          \hline
          $8$ \\ 
          \hline
          $16$\\
          \hline
           $h$\\
          \hline
        \end{tabular}

        \begin{tabular}{|c|c|c|c|c|c|c|c|}
          \hline
          \multicolumn{4}{|c|}{$\nu=-2$} \\ 
          \hline 
          .25 & .33 & .57 & .88\\
          \hline 
          .24 &  .32  & .55 & .63\\
          \hline
          .24 & .31 & .50 & .59 \\
          \hline
          * & .29 & .47& .57 \\
          \hline 
          1/32 & 1/64 & 1/128 & 1/256 \\
          \hline
        \end{tabular}
        \begin{tabular}{|c|c|c|c|c|c|c|c|}
          \hline
          \hline
        \end{tabular}
        \begin{tabular}{|c|c|c|c|c|c|c|c|}
          \hline
          \multicolumn{4}{|c|}{$\nu=-8$} \\
          \hline
          .29 & .47 & .64 & .91\\
          \hline 
          .26 & .57 & .79 & .88\\
          \hline
          .25 & .41 & .71 & .83\\
          \hline
          * &.37 & .68 & .70 \\
          \hline
          1/32 & 1/64 & 1/128 & 1/256\\
          \hline
        \end{tabular} 
      \end{tabular} 
    \end{center}

    \begin{center}\small

       \begin{tabular}{lll}
        \begin{tabular}{|c|c|c|c|c|c|c|c|}
          \hline
          $\tau$\\
          \hline
          $1$\\
          \hline
          $4$\\
          \hline
          $8$ \\ 
          \hline
          $16$\\
          \hline
           $N$\\
          \hline
        \end{tabular}

        \begin{tabular}{|c|c|c|c|c|c|c|c|}
          \hline
          \multicolumn{4}{|c|}{$\nu=-2$} \\ 
          \hline 
          .20 & .21 & .23 & .23\\
          \hline 
          .20 &  .22  & .22 & .22\\
          \hline
          .20 & .21 & .22 & .22 \\
          \hline
          * & .18 & .20& .20 \\
          \hline 
          1/32 & 1/64 & 1/128 & 1/256 \\
          \hline
        \end{tabular}
        \begin{tabular}{|c|c|c|c|c|c|c|c|}
          \hline
          \hline
        \end{tabular}
        \begin{tabular}{|c|c|c|c|c|c|c|c|}
          \hline
          \multicolumn{4}{|c|}{$\nu=-8$} \\
          \hline
          .20 & .21 & .23 & .23\\
          \hline 
          .21 & .22 & .22 & .22\\
          \hline
          .21 & .22 & .23 & .23\\
          \hline
          * &.18 & .20 & .20 \\
          \hline
          1/32 & 1/64 & 1/128 & 1/256\\
          \hline
        \end{tabular} 
      \end{tabular} 
 \caption[$V^{2}$-BootAMG ($k_{r}=8, \eta=4$) for FD Laplace]{
   Asymptotic convergence rate estimates, $\rho$, of the LSR based $V(1,1)$ BAMG solver applied to the FE discretization of the diffusion problem \eqref{eq:jumps}.  The solver is constructed using a $W^2$-cycle BAMG setup with $\eta=4$ and $k_{r}=k_e=8$.  For the results in the bottom tables, an additional intermediate adaptive step is applied to the TV with lowest energy on the finest level only before the second $W$-cycle setup.
\label{tab:aBAMG}}
    \end{center}
  \end{table}

 As the results illustrate, the $W^2$-cycle setup does not yield a scalable solver
 for this problem.\footnote{{In Table 6, * indicates that the fine grid
does not resolve the variations in the permeability field.}}  However, adding an adaptive step to the algorithm before the second $W$-cycle
 does result in an optimal order method.  We note that using additional $W$-cycles in the setup, i.e., a $W^m$-cycle
 with $m > 2$, yields similar results as obtained with this additional adaptive step, but that these 
 additional cycles require significantly more work than does applying the adaptive step a few times only on the finest level, as implemented in the reported experiments.  Generally, using this simple modification reduces the number of bootstrap setup cycles needed to obtain an efficient solver by at least one and in most
 cases two or more for this problem, assuming the number of relaxation steps in the setup cycles is kept fixed.   {We note, in addition, that in~\cite{MBrezina_2005} adaptive AMG was  successfully applied to the same test problem and that an adaptive approach may be preferable for this problem since only a single smooth-error component is needed to computed effective AMG interpolation.}
 
The use of adaptive cycles in the bootstrap setup process has also been considered in various 
 contexts, including the work in~\cite{BAMG_Schwinger} where a BAMG solver was designed for the non-Hermitian Dirac equation of lattice QCD and the work in~\cite{brannick_markov_2011,BAMG_Markov_2}, where a BAMG method for computing state vectors of Markov chains was developed.  These works have shown that
combining the least squares interpolation with the bootstrap setup for computing 
the TVs leads to an efficient and integrated self learning process for defining algebraic multigrid interpolation.
While this research has resulted in an increased understanding of these techniques, a systematic study of how to combine them in an optimal and robust way is still needed.  
 
\section{Future directions and open problems}\label{sec:FD}

Research in the design and analysis of the individual components of the BAMG solver 
has led to efficient and robust multilevel algorithms for a variety of problems, 
including indefinite Helmholtz and Schr\"odinger equations \cite{IL08, IL11}, non-Hermitian PDE systems such as the Dirac equation in lattice QCD~\cite{BAMG_Schwinger}, anisotropic diffusion problems with highly oscillatory and strongly anisotropic 
diffusion coefficient~\cite{BAMG_Aniso}, and the Graph Laplacian problem~\cite{LAMG_Report}.  More generally, this work has led to an increased understanding of the techniques of 
compatible relaxation, algebraic distances, least squares interpolation, 
and the bootstrap process for algebraically computing representatives of
smooth error, the test vectors.   These works have demonstrated the potential
of the BAMG methodology to deliver an efficient and robust algebraic multilevel solver
for a large class of sparse matrix equations.  

The basic components of the BAMG process have been studied and are (mostly) understood
for HPD problems .
The main issue to consider in developing BAMG as a fast solver for general HPD problems concerns the integration of the BAMG components in order to obtain an overall algorithm for use in practice.  Specifically, a systematic study of the BAMG processes for a carefully chosen and diverse set of test problems (that includes
systems of PDEs) is needed.  As noted in the previous sections, 
the most pressing research directions as the authors currently see them are 
\begin{itemize}
\item An efficient  mechanism for  updating coarse variable set $C$ based on a local  slowness of CR relaxation, see Section \ref{sec:CR};  
\item A local procedure for finding weights $\omega_\kappa$ in (\ref{eq:LSfuncrowi}) -- to prevent an emphasize on a globally smooth test vector, that will otherwise be  heavily weighted everywhere,  in the regions  where it has large residuals; such approach will be a generalization of the LSR scheme which is now used concurrently with  global weights, see Section \ref{sec:LS};
\item A balance between work and accuracy when computing  TVs, see Section \ref{sec:BS}; 
\item A procedure  for incorporating all available information about a given problem, in particular, its  geometric information; 
\item An integration of all BAMG components to obtain an easy to tune and robust multigrid solver.
\end{itemize}

Although the basic components of BAMG can in principle also 
be applied to non-Hermitian problems, there are still several basic research questions in this 
setting that require further investigation.  One main unresolved issue that arises in a Galerkin-based algebraic 
multigrid concerns characterizing the algebraically smooth errors that need to be approximated
by restriction and interpolation.  The typical approach has been to compute approximations to singular values with small singular vectors, {e.g., in~\cite{GESA} such a strategy was developed for conviction diffusion problems.}  These approaches are based on the observation that the smallest 
singular value is less than the smallest in magnitude eigenvalue. 
However, for many problems, the smallest eigenvalues of a non-Hermitian matrix are much smaller than the second smallest singular value.  Thus, in general the near kernel of the system matrix can not be characterized
by the singular vectors (or eigenvectors) with small singular values, unless the matrix is normal.

The basic approach in BAMG is to iterate on the homogenous system using the AMG smoother and
the emerging multilevel structure in order to expose the smooth errors that need
to be treated on coarse levels.  Given an appropriate relaxation scheme, 
e.g., the Kaczmarz iteration, which provides a general approach for 
non-Hermitian problems, the subspace of smooth error is in
the most general sense related to the field of values 
\begin{equation*}
  \mathcal F =\left\{ \frac{\langle Ax,x\rangle}{\langle x,x \rangle} \bigg|  x \in \mathbb{C}^n, x \neq 0\right\}.
\end{equation*}   
We note that the field of values may contain the origin even if all eigenvalues of the matrix
have a positive real part.  Thus, it is possible that the Galerkin coarse 
level operator can have a zero eigenvalue, even when the fine-level matrix
does not.  
We have observed in practice that in such
cases, the resulting Galerkin method can 
diverge.
 Here, an underestimate of the small singular values 
or eigenvalues on some grid can lead to an over-correction 
from the coarse grid and divergence of the 
resulting AMG solver.  
Thus, although the BAMG process provides techniques for constructing
restriction $R$ and interpolation $P$, each in their
own setup scheme, these processes must be carefully tied together locally
in order to obtain a suitable Galerkin coarse-grid matrix.

The challenges presented by indefinite problems, e.g., the Helmholtz and Schr\"odinger operators, are again related to the local character of their low eigenspace, which consists of many eigenfunctions that are locally different from one another.  Generally, we refer to such problems as non-elliptic systems.
The key issues is that such eigenfunctions cannot be accurately interpolated by a single prolongation operator (\ref{eq:LSfuncrowi}) with a small caliber, starting with some sufficiently fine scales. There and on all coarser scales several prolongation operators have to be employed, resulting in multiple-Galerkin approach with many coarse-grid operators. To guarantee optimal computational costs, all prolongation operators have to be locally orthogonal, and  achieving that  generates  additional computational difficulties. To this date, the BAMG-like approach has been implemented only for the one-dimensional Helmholtz \cite{IL08} and Schr\"odinger \cite{IL11} problems, which require only {\sl two} prolongation operators (instead of one).   In higher dimensions, the number of prolongation operators and coarse spaces increases, but the ideas carry over in a similar way from 1d.
Two important future research topics are thus to continue to investigate developing BAMG solvers for these and other non-elliptic problems and to determine how far these techniques can be extended as fast solvers for general sparse matrix equations.


\bibliographystyle{plain}
\bibliography{bamg_overview}

\providecommand{\noopsort}[1]{}
\begin{thebibliography}{10}

\bibitem{BAMG1}
A.~Brandt.
\newblock Multi-level adaptive solutions to boundary value problems.
\newblock {\em Math. Comp.}, 31:333--390, 1977.

\bibitem{B83}
A.~Brandt.
\newblock Algebraic multigrid theory: The symmetric case.
\newblock {\em Appl. Math. Comput.}, 19(1-4):23--56, 1986.

\bibitem{etna2000}
A.~Brandt.
\newblock General highly accurate algebraic coarsening.
\newblock {\em Elect. Trans. Numer. Anal.}, 10:1--20, 2000.

\bibitem{B00}
A.~Brandt.
\newblock Multiscale scientific computation: Review 2001.
\newblock In {\em Multiscale and Multiresolution Methods}, pages 1--96.
  Springer Verlag, 2001.

\bibitem{SU}
A.~Brandt.
\newblock Principles of systematic upscaling.
\newblock In {\em Bridging the Scales in Science and Engineering}, pages
  193--215. Oxford University Press, 2010.

\bibitem{brannick_markov_2011}
A.~Brandt, J.~Brannick, M.~Bolten, A.~Frommer, K.~Kahl, and I.~Livshits.
\newblock Bootstrap {AMG} for {M}arkov chains.
\newblock {\em SIAM J. Sci. Comp.}, 33:3425--3446, 2011.

\bibitem{BAMG2010}
A.~Brandt, J.~Brannick, K.~Kahl, and I.~Livshits.
\newblock Bootstrap {AMG}.
\newblock {\em SIAM J. Sci. Comput.}, 33:612--632, 2011.

\bibitem{BAMG_Aniso}
A.~Brandt, J.~Brannick, K.~Kahl, and I.~Livshits.
\newblock Algebraic distance as a measure of strength of connection in {AMG}.
\newblock {\em Electronic Transactions in Numerical Analysis}, 2013.
\newblock Submitted: April 2013. Also available as arXiv:1106.5990 [math.NA].

\bibitem{BMR83}
A.~Brandt, S.~McCormick, and J.~W. Ruge.
\newblock Algebraic multigrid ({AMG}) for automatic multigrid solution with
  application to geodetic computations.
\newblock Technical report, Colorado State University, Fort Collins, Colorado,
  1983.

\bibitem{BMR84}
A.~Brandt, S.~F. McCormick, and J.~W. Ruge.
\newblock Algebraic multigrid ({AMG}) for sparse matrix equations.
\newblock In D.~J. Evans, editor, {\em Sparsity and Its Applications}.
  Cambridge University Press, Cambridge, 1984.

\bibitem{RMG}
A.~Brandt and D.~Ron.
\newblock Renormalization multigrid ({RMG}): Statistically optimal
  renormalization group flow and coarse-to-fine monte carlo acceleration.
\newblock {\em Journal of Statistical Physics}, 102:231--257, 2001.
\newblock 10.1023/A:1026520927784. First appeared as Renormalization Multigrid
  ({RMG}): Statistically Optimal Renormalization Group Flow and Coarse-to-Fine
  {Monte-Carlo} Acceleration, Technical Report GMC-11, Weizmann Institute of
  Science, 1999.

\bibitem{JBrannick_2005a}
J.~Brannick.
\newblock {\em Adaptive algebraic Multigrid coarsening strategies}.
\newblock PhD thesis, University of Colorado at Boulder, 2005.

\bibitem{brannick_local_stab_2011}
J.~Brannick, Y.~Chen, and L.~Zikatanov.
\newblock An algebraic multilevel method for anisotropic elliptic equations
  based on subgraph matching.
\newblock {\em Numer. Linear Algebra Appl.}, To appear, accepted for
  publication November 17, 2011.

\bibitem{JBrannick_RFalgout}
J.~Brannick and R.~Falgout.
\newblock Compatible relaxation and coarsening in algebraic multigrid.
\newblock {\em SIAM J. Sci. Comput.}, 32(3):1393--1416, 2010.

\bibitem{BAMG_Schwinger}
J.~Brannick and K.~Kahl.
\newblock Bootstrap algebraic multigrid for the 2d wilson-dirac system.
\newblock {\em SIAM Journal of Scientific Computing}.
\newblock Submitted (August 28, 2013). Also available on aRxiV:1308.5992
  [math.NA].

\bibitem{BAMG_Markov_2}
J.~Brannick, K.~Kahl, and S.~Sokolovic.
\newblock {\em Journal of Applied Numerical Mathematics}, Submitted January 19,
  2014.

\bibitem{Brannick_Trace_06}
J.~Brannick and L.~Zikatanov.
\newblock Algebraic multigrid methods based on compatible relaxation and energy
  minimization.
\newblock In O.~B. Widlund and D.~E. Keyes, editors, {\em Lecture Notes in
  Computational Science and Engineering}, volume~55, pages 15--26. Springer,
  2007.

\bibitem{MBrezina_RFalgout_SMacLachlan_TManteuffel_SMcCormick_JRuge_2003}
M.~Brezina, R.~Falgout, S.~MacLachlan, T.~Manteuffel, S.~McCormick, and
  J.~Ruge.
\newblock Adaptive smoothed aggregation ($\alpha${SA}).
\newblock {\em SIAM J. Sci. Comput.}, 25(6):1896--1920, 2004.

\bibitem{MBrezina_2005}
M.~Brezina, R.~Falgout, S.~MacLachlan, T.~Manteuffel, S.~McCormick, and
  J.~Ruge.
\newblock Adaptive amg (a{AMG}).
\newblock {\em SIAM J. Sci. Comput.}, 26:1261--1286, 2005.

\bibitem{Brezina:2012:RCB}
M.~Brezina, C.~Ketelsen, T.~Manteuffel, S.~McCormick, M.~Park, and J.~Ruge.
\newblock Relaxation-corrected bootstrap algebraic multigrid ({$r$BAMG)}.
\newblock {\em Journal of Numerical Linear Algebra with Applications},
  19(2):178--193, March 2012.

\bibitem{GESA}
M.~Brezina, T.~Manteuffel, S.~McCormick, J.~Ruge, G.~Sanders, and
  P.~Vassilevski.
\newblock A generalized eigensolver based on smoothed aggregation (ges-sa) for
  initializing smoothed aggregation (sa) multigrid.
\newblock {\em Numer. Linear Algebra Appl.}, 15:249--269, 2008.

\bibitem{rhoAMGe}
T.~Chartier, R.~D. Falgout, V.~E. Henson, J.~E.~Jones a~nd T.~Manteuffel, S.~F.
  McCormick, J.~W. Ruge, and P.~S. Vassilevski.
\newblock Spectral {AMG}e ($\rho${AMG}e).
\newblock {\em SIAM J. Sci. Comput.}, 25:1--26, 2003.

\bibitem{PanayotRob_2003}
R.~Falgout and P.~Vassilevski.
\newblock On generalizing the amg framework.
\newblock {\em SIAM J. Numer. Anal.}, 42(4):1669--1693, 2004.

\bibitem{antrg}
N.~Halko, P.~Martinsson, and J.~Tropp.
\newblock Finding structure with randomness: Probabilistic algorithms for
  constructing approximate matrix decompositions.
\newblock {\em SIAM Review}, 53:217Ð288, 2011.

\bibitem{HCR}
O.~E. Livne.
\newblock Coarsening by compatible relaxation.
\newblock {\em Num. Lin. Alg. Appl.}, 11:205--227, 2004.

\bibitem{LAMG_Report}
O.E. Livne and A.~Brandt.
\newblock Lean algebraic multigrid ({LAMG}): Fast graph {L}aplacian linear
  solver.
\newblock {\em SIAM Journal of Scientific Computing}, 34:B499--B522, 2012.

\bibitem{IL08}
I.~Livshits.
\newblock One-dimensional algorithm for finding eigenbasis of the
  schr{\"o}dinger operator.
\newblock {\em SIAM J. Sci. Comput.}, 30(1):416--440, 2008.

\bibitem{IL11}
I.~Livshits.
\newblock The least squares amg solver for the one-dimensional helmholtz
  operator.
\newblock {\em Computing and Visualization in Science}, (14):17--25, 2011.

\bibitem{MacRed}
S.~MacLachlan, T.~Manteuffel, and S.~McCormick.
\newblock Adaptive reduction-based amg.
\newblock {\em Numer. Linear Algebra Appl.}, 13:599Ð620, 2006.

\bibitem{SMacLYS07}
S.~MacLachlan and Y.~Saad.
\newblock A greedy strategy for coarse-grid selection.
\newblock {\em SIAM J. Sci. Comput.}, 29(5):1825Ð1853, 2007.

\bibitem{iBAMG}
T.~Manteuffel, S.~McCormick, M.~Park, and J.~Ruge.
\newblock Operator-based interpolation for bootstrap algebraic multigrid.
\newblock {\em Numer. Linear Algebra Appl.}, 17(2-3):519--537, 2010.

\bibitem{Adap_Smooth}
Bobby Philip and Timothy~P. Chartier.
\newblock Adaptive algebraic smoothers.
\newblock {\em J. Computational Applied Mathematics}, 236:2277--2297, 2012.

\bibitem{RonSB11}
D.~Ron, I.~Safro, and A.~Brandt.
\newblock Relaxation-based coarsening and multiscale graph organization.
\newblock {\em Multiscale Modeling and Simulation}, 9:407--423, 2011.

\bibitem{RS87}
J.~W. Ruge and K.~St{\"u}ben.
\newblock Algebraic multigrid ({AMG}).
\newblock In S.~F. McCormick, editor, {\em Multigrid Methods}, volume~3 of {\em
  Frontiers in Applied Mathematics}, pages 73--130. SIAM, Philadelphia, PA,
  1987.

\bibitem{SRB09}
I.~Safro, D.~Ron, and A.~Brandt.
\newblock Multilevel algorithms for linear ordering problems.
\newblock {\em J. Exp. Algorithmics}, 13:1.4--1.20, 2009.

\bibitem{Jacob}
J.~B. Schroder.
\newblock Smoothed aggregation solvers for anisotropic diffusion.
\newblock {\em Numer. Linear Algebra Appl.}, 19:296--312, 2012.

\bibitem{JWilkinson_1965}
J.~Wilkinson.
\newblock {\em The Algebraic eigenvalue problem}.
\newblock Clarendon Press, Oxford, 1965.

\end{thebibliography}

\end{document}